\documentclass[a4paper]{amsart}
\usepackage{amssymb}
\usepackage{verbatim}
\usepackage{pb-diagram}
\usepackage{epic,eepic}
\usepackage[dvips]{graphicx}
\theoremstyle{plain}
  \newtheorem{thm}{Theorem}[section]
  
  \newtheorem{cor}[thm]{Corollary}
  
  \newtheorem{conj}[thm]{Conjecture}

  \newtheorem*{obs*}{Observation}
\theoremstyle{definition}
  
  \newtheorem{ex}[thm]{Example}

\theoremstyle{remark}
  \newtheorem{rem}[thm]{Remark}
  
  \newtheorem*{ack}{Acknowledgments}
\newcommand{\Z}{\mathbb{Z}}
\newcommand{\C}{\mathbb{C}}
\newcommand{\R}{\mathbb{R}}

\newcommand{\Vol}{\operatorname{Vol}}
\newcommand{\CS}{\operatorname{CS}}
\newcommand{\cs}{\operatorname{cs}}
\newcommand{\Li}{\operatorname{Li}}

\newcommand{\arccosh}{\operatorname{arccosh}}

\newcommand{\Hom}{\operatorname{Hom}}
\newcommand{\length}{\operatorname{length}}
\newcommand{\torsion}{\operatorname{torsion}}
\newcommand{\tr}{\operatorname{tr}}

\renewcommand{\Re}{\operatorname{Re}}
\renewcommand{\Im}{\operatorname{Im}}
\numberwithin{equation}{section}
\begin{document}
\title
{An introduction to the volume conjecture and its generalizations}
\author{Hitoshi Murakami}
\address{
Department of Mathematics,
Tokyo Institute of Technology,
Oh-okayama, Meguro, Tokyo 152-8551, Japan
}
\email{starshea@tky3.3web.ne.jp}
\date{\today}
\begin{abstract}
In this paper we give an introduction to the volume conjecture and its generalizations.
Especially we discuss relations of the asymptotic behaviors of the colored Jones polynomials of a knot with different parameters to representations of the fundamental group of the knot complement at the special linear group over complex numbers by taking the figure-eight knot and torus knots as examples.
\end{abstract}
\thanks{This work is supported by Grant-in-Aid for Exploratory Research (18654009) and Travel Grant for Academic Meetings, Japan Society for the Promotion of Science.}
%%%%%%%%%%%%%%%%%%%%%%%%%%%%%%%%%%%%%%%%%%%%%%%%%%%%%%%%%%%%%%%%%%%%%%%%%
\maketitle
%\tableofcontents
%%%%%%%%%%%%%%%%%%%%%%%%%%%%%%%%%%%%%%%%%%%%%%%%%%%%%%%%%%%%%%%%%%%%%%%%%%
After V.~Jones' discovery of his celebrated polynomial invariant $V(K;t)$ in 1985 \cite{Jones:BULAM385}, Quantum Topology has been attracting many researchers; not only mathematicians but physicists.
The Jones polynomial was generalized to two kinds of two-variable polynomials, the HOMFLYpt polynomial \cite{HOMFLY:BULAM385,Przytycki/Traczyk:KOBJM88} and the Kauffman polynomial \cite{Kauffman:Knots} (see also \cite{Ho:ABSAM1985,Brandt/Lickorish/Millett:INVEM1986} for another one-variable specialization).
It turned out that these polynomial invariants are related to quantum groups introduced by V.~Drinfel$'$d and M.~Jimbo (see for example \cite{Kassel:quantum_groups, Turaev:quantum}) and their representations.
For example the Jones polynomial comes from the quantum group $U_q(sl_2(\C))$ and its two-dimensional representation.
We can also define the quantum invariant associated with a quantum group and its representation.
\par
If we replace the quantum parameter $q$ of a quantum invariant ($t$ in $V(K;t)$) with $e^h$ we obtain a formal power series in the formal parameter $h$.
Fixing a degree $d$ of the parameter $h$, all the degree $d$ coefficients of quantum invariants share a finiteness property.
By using this property, one can define a notion of finite type invariant \cite{Birman/Lin:INVEM93,BarNatan:TOPOL95}.
(See \cite{Vassiliev:1990} for V.~Vassiliev's original idea.)
\par
Via the Kontsevich integral \cite{Kontsevich:1993} (see also \cite{Le/Murakami:COMPO11996}) we can recover a quantum invariant from the corresponding `classical' data.
(Note that a quantum group is a deformation of a `classical' Lie algebra.)
So for example one only needs to know the Lie algebra $sl_2(\C)$ (easy!) and its fundamental two-dimensional representation (very easy!) to define the Jones polynomial, provided that one knows the Kontsevich integral (unfortunately, this is very difficult).
\par
In the end of the 20th century, M.~Khovanov introduced yet another insight to Quantum Topology.
He categorified the Jones polynomial and defined a homology for a knot such that its graded Euler characteristic coincides with the Jones polynomial \cite{Khovanov:DUKMJ2000}.
See \cite{Khovanov/Rozansky:2004} for a generalization to the HOMFLYpt polynomial.
\par
Now we are in the 21st century.
\par
In 2001, J.~Murakami and the author proposed the Volume Conjecture \cite{Murakami/Murakami:ACTAM12001} to relate a series of quantum invariants, the $N$-colored Jones polynomials of a knot, to the volume of the knot complement, generalizing R.~Kashaev's conjecture \cite{Kashaev:LETMP97}.
The aim of this article is to introduce the conjecture and some of its generalizations, emphasizing a relation of the asymptotic behavior of the series of the colored Jones polynomials of a knot to representations of the fundamental group of the knot complement.
\par
In Section~1 we prepare some fundamental facts about the colored Jones polynomial, the character variety, the volume and the Chern--Simons invariant.
Section~2 contains our conjectures, and in Sections~3 and 4 we give supporting evidence for the conjectures taking the figure-eight knot and torus knots as examples.
%%%%%%%%%%%%%%%%%%%%%%%%%%%%%%%%%%%%%%%%%%%%%%%%%%%%%%%%%%%%%%%%%%%%
\begin{ack}
The author would like to thank the Institute of Mathematics, Hanoi for its hospitality.
He also thanks the organizers of the International Conference on Quantum Topology from 6th to 12th August, 2007.
\par
He is also grateful to Akiko Furuya for her assistance and encouragement about Travel Grant for Academic Meetings of JSPS.
\end{ack}
%%%%%%%%%%%%%%%%%%%%%%%%%%%%%%%%%%%%%%%%%%%%%%%%%%%%%%%%%%%%%%%%%%%%%%%%%%
\section{Preliminaries}
In this section we review the definition of the colored Jones polynomial, the character variety, and the volume and the Chern--Simons invariant.
\subsection{Colored Jones polynomial}
\begin{figure}[h]
\includegraphics[scale=0.3]{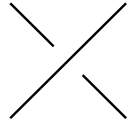}
\caption{A crossing of a knot diagram.}
\label{fig:crossing}
\end{figure}

Let $K$ be an oriented knot in the three-sphere $S^3$, and $D$ its diagram.
We assume that $D$ is the image of a projection $\R^3\cong S^3\setminus\{\infty\}\to\R^2$ of a circle $S^1$ together with a finite number of double point singularities as in Figure~\ref{fig:crossing}.
%%%%%%%%%%%%%%%%%%
%%%%%%%%%%%%%%%%%
We call these double points with over/under data crossings.
\par
Now we associate the knot diagram to a Laurent polynomial in $A$ following L.~Kauffman \cite{Kauffman:TOPOL87}.
\par
First of all we forget the orientation of $D$.
\par
Then replace a crossing with a linear combination of two nullified crossings:
\begin{equation*}
  \raisebox{-6.5mm}{\includegraphics[scale=0.3]{crossing_a.eps}}
  =
  A
  \raisebox{-6.5mm}{\includegraphics[scale=0.3]{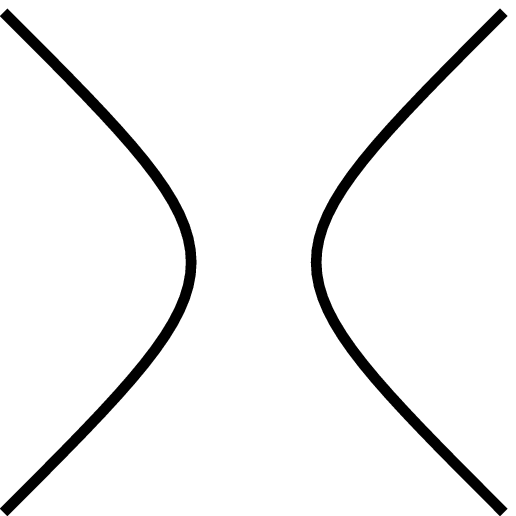}}
  +
  A^{-1}
  \raisebox{-6.5mm}{\includegraphics[scale=0.3]{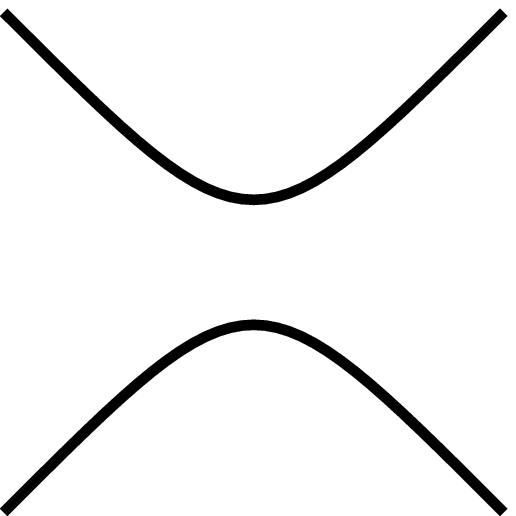}}.
\end{equation*}
Now we have a linear combination of two new diagrams, which may not be knot diagrams but link diagrams.
\par
We choose another crossing and replace it with a linear combination of two new diagrams.
Now we have a linear combination of four diagrams.
\par
Continue these processes until we get a linear combination of diagrams with no crossings, whose coefficients are some powers of $A$.
\par
If we replace each resulting diagram with $(-A^2-A^{-2})^{\nu-1}$ where $\nu$ is the number of components of the diagram, we get a Laurent polynomial $\langle{D}\rangle\in\Z[A,A^{-1}]$, which is called the Kauffman bracket.
\par
\begin{figure}[h]
\raisebox{-6.5mm}{\includegraphics[scale=0.3]{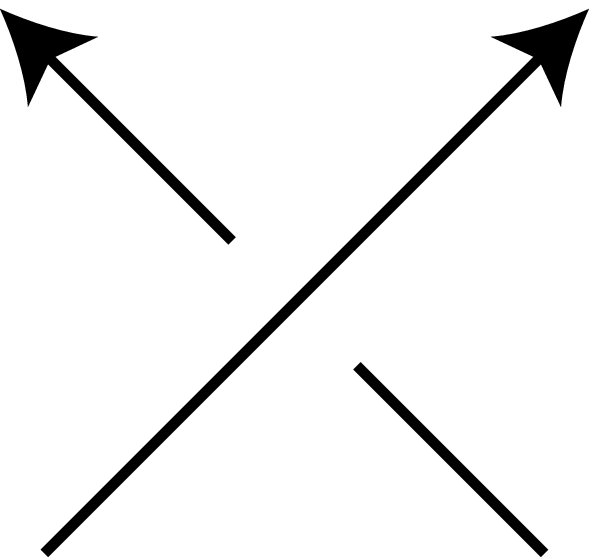}}\quad:$+1$
\quad,\qquad
\raisebox{-6.5mm}{\includegraphics[scale=0.3]{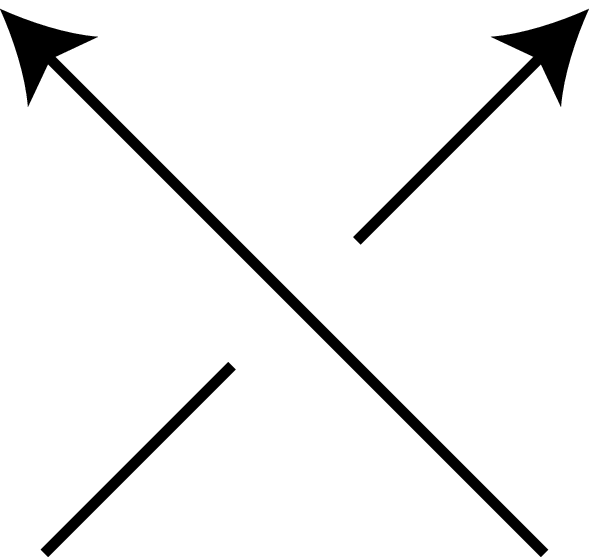}}\quad:$-1$
\caption{A positive crossing (left) and a negative crossing (right).}
\label{fig:sign_crossing}
\end{figure}

Now we recall the orientation of $D$.
Let $w(D)$ be the sum of the signs of the crossings in $D$, where the sign is defined as in Figure~\ref{fig:sign_crossing}.
%%%%%%%%%%%%%%%%%%
%%%%%%%%%%%%%%%%%%
\par
Lastly we define $V(K;t):=\left(-A^3\right)^{-w(D)}\langle{D}\rangle\Big|_{A:=t^{-1/4}}$.
It is known that this coincides with (a version of) the Jones polynomial \cite{Jones:BULAM385} since
\begin{equation*}
\begin{split}
  &t^{-1}V
  \left(
    \raisebox{-2mm}{\includegraphics[scale=0.1]{crossing_p_a.eps}}\,;t
  \right)
  -
  tV
  \left(
    \raisebox{-2mm}{\includegraphics[scale=0.1]{crossing_n_a.eps}}\,;t
  \right)
  \\
  =&
  t^{-1}\times
  \left(-t^{-3/4}\right)^{-w
    \left(
      \raisebox{-0.7mm}{\includegraphics[scale=0.05]{crossing_p_a.eps}}
    \right)}
  \left\langle
    \raisebox{-1.8mm}{\includegraphics[scale=0.1]{crossing_a.eps}}
  \right\rangle
  -
  t\times
  \left(-t^{-3/4}\right)^{-w
    \left(
      \raisebox{-0.7mm}{\includegraphics[scale=0.05]{crossing_n_a.eps}}
    \right)}
  \left\langle
    \raisebox{-1.8mm}{\includegraphics[scale=0.1]{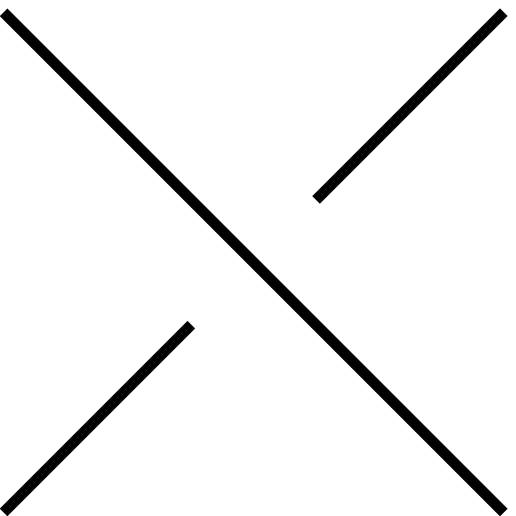}}
  \right\rangle
  \\
  =&
  t^{-1}\times
  \left(-t^{-3/4}\right)^{-w
    \left(
      \raisebox{-0.7mm}{\includegraphics[scale=0.05]{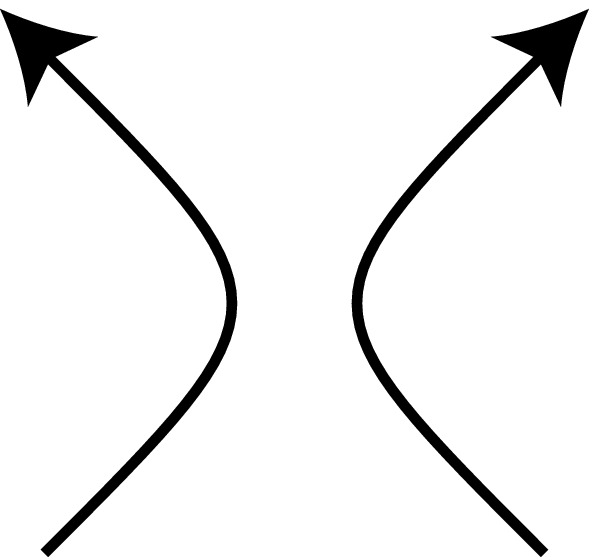}}
    \right)-1}
  \left\{
    t^{-1/4}
    \left\langle
      \raisebox{-1.8mm}{\includegraphics[scale=0.1]{crossing_v_a.eps}}
    \right\rangle
    +t^{1/4}
    \left\langle
      \raisebox{-1.8mm}{\includegraphics[scale=0.1]{crossing_h_a.eps}}
    \right\rangle
  \right\}
  \\
  &-
  t\times
  \left(-t^{-3/4}\right)^{-w
    \left(
      \raisebox{-0.7mm}{\includegraphics[scale=0.05]{crossing_0_a.eps}}
    \right)+1}
  \left\{
    t^{-1/4}
    \left\langle
      \raisebox{-1.8mm}{\includegraphics[scale=0.1]{crossing_h_a.eps}}
    \right\rangle
    +t^{1/4}
    \left\langle
      \raisebox{-1.8mm}{\includegraphics[scale=0.1]{crossing_v_a.eps}}
    \right\rangle
  \right\}
  \\
  =&
  \left(-t^{-3/4}\right)^{-w
    \left(
      \raisebox{-0.7mm}{\includegraphics[scale=0.05]{crossing_0_a.eps}}
    \right)}
  \left(
    t^{1/2}-t^{-1/2}
  \right)
  \left\langle
    \raisebox{-1.8mm}{\includegraphics[scale=0.1]{crossing_v_a.eps}}
  \right\rangle
  \\
  =&
  \left(
    t^{1/2}-t^{-1/2}
  \right)V
  \left(
    \raisebox{-2mm}{\includegraphics[scale=0.1]{crossing_0_a.eps}}\,;t
  \right).
\end{split}
\end{equation*}
\par
In this paper we use another version $J_2(K;q)$ that satisfies
\begin{equation*}
  qJ_2
  \left(
    \raisebox{-2mm}{\includegraphics[scale=0.1]{crossing_p_a.eps}}\,;q
  \right)
  -
  q^{-1}J_2
  \left(
    \raisebox{-2mm}{\includegraphics[scale=0.1]{crossing_n_a.eps}}\,;q
  \right)
  =
  \left(
    q^{1/2}-q^{-1/2}
  \right)J_2
  \left(
    \raisebox{-2mm}{\includegraphics[scale=0.1]{crossing_0_a.eps}}\,;q
  \right)
\end{equation*}
with the normalization $J_2(\text{unknot};q)=1$.
\par
It is well known that $J_2$ corresponds to the two-dimensional irreducible representation of the Lie algebra $sl(2;\C)$.
The invariant corresponding to the $N$-dimensional irreducible representation is called the $N$-colored Jones polynomial and denoted by $J_N(K;q)$.
It can be also defined as a linear combination of the ($2$-colored) Jones polynomials of $(N-1)$-parallels or less of the original knot.
\par
\begin{figure}[h]
\includegraphics[scale=0.5]{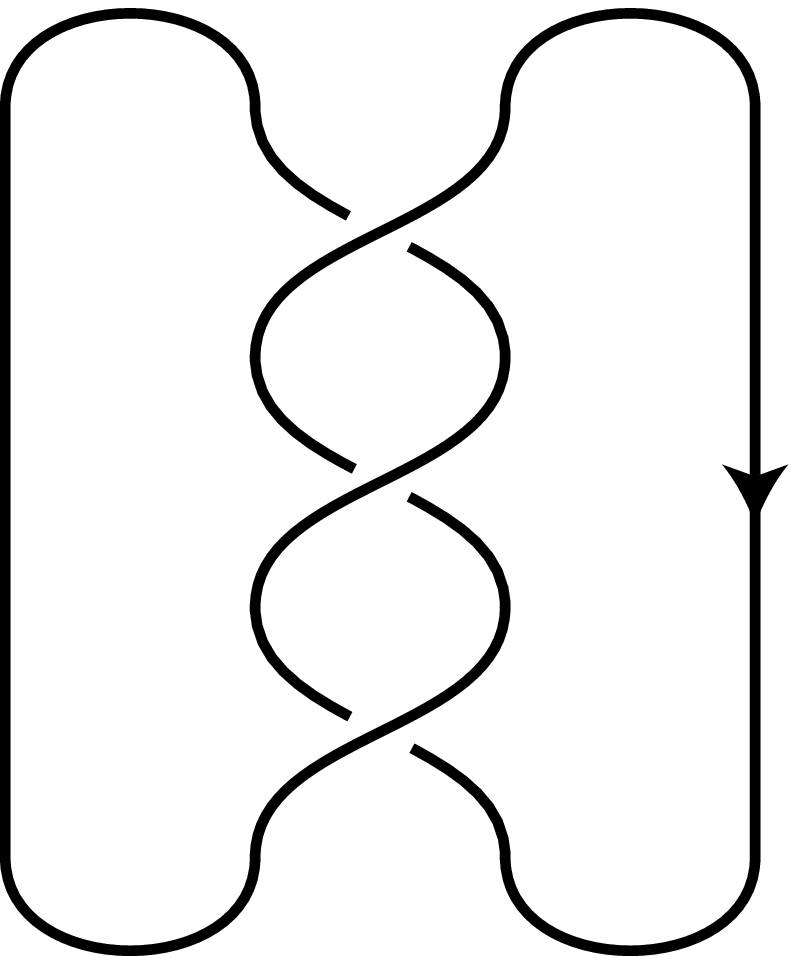}
\caption{Trefoil.}
\label{fig:trefoil}
\end{figure}
\begin{figure}[h]
\includegraphics[scale=0.5]{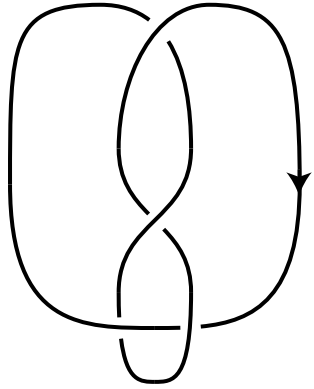}
\caption{Figure-eight knot}
\label{fig:fig8}
\end{figure}

For example, K.~Habiro \cite{Habiro:SURIK00} and T.~L{\^e} \cite{Le:TOPOA2003} calculated the $N$-colored Jones polynomial for the trefoil knot $T$ (Figure~\ref{fig:trefoil}) and for the figure-eight knot $E$ (Figure~\ref{fig:fig8}) as follows:
\begin{equation*}
  J_N(T;q)
  =
  q^{1-N}
  \sum_{k=0}^{N-1}
  q^{-kN}
  \prod_{j=1}^{k}
  \left(1-q^{j-N}\right)
\end{equation*}
and
\begin{equation*}
  J_N(E;q)
  =
  \sum_{k=0}^{N-1}
  \prod_{j=1}^{k}
  \left(q^{(N-j)/2}-q^{-(N-j)/2}\right)
  \left(q^{(N+j)/2}-q^{-(N+j)/2}\right).
\end{equation*}
For other formulas of $J_N(T;q)$, see \cite{Morton:MATPC95,Masbaum:ALGGT12003}.
%%%%%%%%%%%%%
%%%%%%%%%%%%%
%%%%%%%%%%%%%%%%%%%%%%%%%%%%%%%%%%%%%%%%%%
\subsection{$SL(2;\C)$-Character variety of a knot group}\label{subsec:character_variety}
For a knot in $S^3$, let us consider representations of $\pi_1(S^3\setminus{K})$ at $SL(2;\C)$, where $SL(2;\C)$ is the set of all $n\times n$ complex matrices with determinant one.
The set of all characters of representations is denoted by $X(S^3\setminus{K})$ or $X(K)$.
We can regard $X(K)$ as the set of equivalence classes where we regard two representations equivalent if they have the same trace.
It is well known that $X(K)$ is an algebraic variety (the set of the zeroes of polynomials) and called the $SL(2;\C)$-character variety of the knot $K$.
(See for example \cite{Cooper/Culler/Gillet/Long/Shalen:INVEM1994}.)
\par
Here are two families of knots whose $SL(2;\C)$-character varieties are well known.
%%%%%%%%%%%%%%%%%%%%%%%%%%%%%%%%
\begin{ex}[Two-bridge knots]
If $K$ is a two-bridge knot (a knot that can be put in $\R^3=S^3\setminus\{\infty\}$ in such a way that it has only two maxima), we can express its $SL(2;\C)$-character variety as follows \cite[3.3.1~Theorem]{Le:RUSSM1991} (see also \cite[Theorem~1]{Riley:QUAJM31984}).
\par
It is well known that $\pi_1(S^3\setminus{K})$ has the following presentation:
\begin{equation*}
  \langle
    x,y\mid \omega x=y\omega
  \rangle,
\end{equation*}
where $\omega$ is a word in $x$ and $y$ of the form $x^{\pm1}y^{\pm1}\dots x^{\pm1}y^{\pm1}$.
Let $\rho$ be a representation of $\pi_1(S^3\setminus{K})$ at $SL(2;\C)$.
Then for any word $z$ in $x$ and $y$ we can express $\tr\bigl(\rho(z)\bigr)$ as a polynomial in $\xi:=\tr\bigl(\rho(x)\bigr)$ and $\eta:=\tr\bigl(\rho(xy)\bigr)$ by using the following formulas \cite{Whittemore:PROAM1973}:
\begin{equation}\label{eq:trace}
\begin{split}
  \tr(AB)&=\tr(A)\tr(B)-\tr(A^{-1}B), \\
  \tr(A^{-1})&=\tr(A).
\end{split}
\end{equation}
We denote this polynomial by $P_{z}(\xi,\eta)$.
\par
The $SL(2;\C)$-character variety $X(K)$ is given as follows.
\begin{thm}[Le \cite{Le:RUSSM1991}]\label{thm:Le}
The $SL(2;\C)$-character variety is determined by the polynomial
\begin{equation*}
  (2+\eta-\xi^2)F(\xi,\eta),
\end{equation*}
where
\begin{equation*}
  F(\xi,\eta)
  :=
  \sum_{i=0}^{k}(-1)^{i}P_{\omega^{(i)}}(\xi,\eta)
\end{equation*}
with $\omega^{(i)}$ the word obtained from $\omega$ by deleting the first $i$ letters and the last $i$ letters (we put $P_{\emptyset}(\xi,\eta):=1$).
Moreover the first factor $2+\eta-\xi^2$ determines the abelian part and the second factor $F(\xi,\eta)$ determines the non-abelian part.
\end{thm}
\end{ex}
%%%%%%%%%%%
\begin{ex}[Torus knots]
A knot that can be put on the standard torus is called a torus knot.
Such knots are parametrized by two coprime integers, and up to mirror image we may assume that they are both positive.
We also assume that they are bigger than $1$.
\par
For the torus knot $T(a,b)$ of type $(a,b)$ with positive coprime integers $a$ and $b$ with $a>1$ and $b>1$, the fundamental group $\pi_1(T(a,b))$ also has a presentation with two generators and one relation:
\begin{equation*}
  \langle
    g,h\mid g^a=h^b
  \rangle.
\end{equation*}
Note that the meridian (a loop that goes around the knot in the positive direction) $\mu$ and the longitude (a loop that is parallel to the knot and is null-homologous in the knot complement) $\lambda$ can be expressed as $\mu=g^{-c}h^{d}$ and $\lambda=g^{a}\mu^{ab}$, where we choose $c$ and $d$ so that $ad-bc=1$.
\par
Then the $PSL(2;\C)$-character variety of $\pi_1(T(a,b))$ is given as follows.
\begin{thm}[{\cite[Theorem~2]{Dubois/Kashaev:MATHA2007}}]\label{thm:Dubois_Kashaev}
The components of the non-abelian part of the $SL(2;\C)$-character variety are indexed by $k$ and $l$ with $1\le k \le a-1$, $1\le l \le b-1$ and $k\equiv l\pmod{2}$.
Moreover a representation $\rho_{k,l}$ in the component indexed by $(k,l)$ satisfies
\begin{align*}
  \tr\bigl(\rho_{k,l}(g)\bigr)
  &=
  2\cos\left(\frac{k\pi}{a}\right)
  \intertext{and}
  \tr\bigl(\rho_{k,l}(h)\bigr)
  &=
  2\cos\left(\frac{l\pi}{b}\right).
\end{align*}
\end{thm}
(See also \cite[Theorem~1]{Klassen:TRAAM1991}.)
\end{ex}
%%%%%%%%%%%%%%%%%%%%%%%%%%%%%%%%%%%%%%%%%%%%%%%%%%
\subsection{Volume and the Chern--Simons invariant of a representation}
\label{subsec:V_CS}
For a closed three-manifold $W$ one can define the Chern--Simons function $\cs_{M}$ on the $SL(2;\C)$-character variety of $\pi_1(W)$ as follows.
See for example \cite[Section~2]{Kirk/Klassen:MATHA90} for a nice review.
\par
Let $A$ be a $\mathfrak{sl}(2;\C)$-valued $1$-form on $W$ satisfying $d\,A+A\wedge A=0$.
Then $A$ defines a flat connection on $W\times{SL(2;\C)}$ and so it also induces the holonomy representation $\rho$ of $\pi_1(W)$ at $SL(2;\C)$ up to conjugation.
Then the $SL(2;\C)$ Chern--Simons function $\cs_{W}$ is defined by
\begin{equation*}
  \cs_{W}\bigl([\rho]\bigr)
  :=
  \frac{1}{8\pi^2}
  \int_{W}\tr\left(A\wedge d\,A+\frac{2}{3}A\wedge A\wedge A\right)
  \in\C
  \pmod{\Z},
\end{equation*}
where $[\rho]$ means the conjugacy class.
Note that every representation is induced as the holonomy representation of a flat connection.
\par
P.~Kirk and E.~Klassen \cite{Kirk/Klassen:COMMP93} gave a formula to calculate the $SL(2;\C)$ Chern--Simons invariant if $W$ is given as the union of two three-manifolds whose boundaries are tori.
\par
Let $M$ be an oriented three-manifold with boundary $\partial{M}$ a torus.
Given a representation $\rho$ of $\pi_1(M)$ at $SL(2;\C)$, one can define the Chern--Simons function $\cs_{M}$ as a map from the $SL(2;\C)$-character variety $\tilde{X}(M)$ of $\pi_1(M)$ to a circle bundle $E(\partial{M})$ over the character variety $\tilde{X}(\partial{M})$ of $\pi_1(\partial{M})\cong\Z\oplus\Z$ \cite{Kirk/Klassen:COMMP93}, which is a lift of the restriction map $\tilde{X}(M)\to\tilde{X}(\partial{M})$.
\par
To describe $\cs_{M}$, we consider a map
\begin{equation*}
  p\colon\Hom(\pi_1(\partial{M}),\C)\to\tilde{X}(\partial{M})
\end{equation*}
defined by
\begin{equation*}
  p(\kappa)
  :=
  \left[
    \gamma\mapsto
    \begin{pmatrix}
      \exp(2\pi\sqrt{-1}\kappa(\gamma)) & 0 \\
      0 & \exp(-2\pi\sqrt{-1}\kappa(\gamma))
    \end{pmatrix}
  \right]
\end{equation*}
for $\gamma\in\pi_1(\partial{M})$, where the square brackets mean the equivalence class in $\tilde{X}(\partial{M})$.
Note that $\kappa$ and $\kappa'$ define the same element in $\tilde{X}(\partial{M})$ if and only if $\cos(2\pi\sqrt{-1}\kappa(\gamma))=\cos(2\pi\sqrt{-1}\kappa'(\gamma))$ for any $\gamma\in\pi_1(\partial{M})$.
\par
Now fix a generator system $(\mu,\lambda)$ of $\pi_1(\partial{M})\cong\Z\oplus\Z$, and take its dual basis $(\mu^{\ast},\lambda^{\ast})$ of $\Hom(\pi_1(\partial{M}),\C)$.
If $\kappa=\alpha\mu^{\ast}+\beta\lambda^{\ast}$ and $\kappa'=\alpha'\mu^{\ast}+\beta'\lambda^{\ast}$, then $p(\kappa)=p(\kappa')$ if and only if
\begin{equation*}
  k\alpha+l\beta=\pm(k\alpha'+l\beta')\pmod{\Z}
\end{equation*}
for any integers $k$ and $l$.
This means that $p$ is invariant under the following actions on $\Hom(\pi_1(\partial{M});\C)$:
\begin{align*}
  x\cdot(\alpha,\beta)&:=(\alpha+1,\beta),
  \\
  y\cdot(\alpha,\beta)&:=(\alpha,\beta+1),
  \\
  b\cdot(\alpha,\beta)&:=(-\alpha,-\beta).
\end{align*}
Note that these actions form the group
\begin{equation*}
  G
  :=
  \langle
    x,y,b\mid
    xyx^{-1}y^{-1}=bxbx=byby=b^2=1
  \rangle
\end{equation*}
and that in fact the quotient space $\Hom(\pi_1(\partial{M}),\C)/G$ can be identified with the $SL(2;\C)$-character variety $\tilde{X}(\partial{M})$.
\par
Next let us consider the following actions on $\Hom(\pi_1(\partial{M}),\C)\times\C^{\ast}$.
\begin{align*}
  x\cdot[\alpha,\beta;z]
  &:=
  [\alpha+1,\beta;z\exp(2\pi\sqrt{-1}\beta)],
  \\
  y\cdot[\alpha,\beta;z]
  &:=
  [\alpha,\beta+1;z\exp(-2\pi\sqrt{-1}\alpha)],
  \\
  b\cdot[\alpha,\beta;z]
  &:=
  [-\alpha,-\beta;z].
\end{align*}
We denote the quotient space $(\Hom(\pi_1(\partial{M}),\C)\times\C^{\ast})/G$ by $E(\partial{M})$.
Then $E(\partial{M})$ becomes a $\C^{\ast}$-bundle over $\tilde{X}(\partial{M})$.
\par
The Chern--Simons function $\cs_M$ is a map from $\tilde{X}(M)$ to $E(\partial{M})$ such that the following diagram commutes, where $q$ is the projection $q\colon E(\partial{M})\to\tilde{X}(\partial{M})$ and $i^{\ast}$ is the restriction map.
\begin{equation*}
\begin{diagram}
  \node[2]{E(\partial{M})}\arrow{s,r}{q} \\
  \node{\tilde{X}(M)}\arrow{ne,t}{\cs_M}\arrow{e,t}{i^{\ast}}\node{\tilde{X}(\partial{M})}
\end{diagram}
\end{equation*}
\par
If a closed three-manifold $W$ is given as $M_1\cup M_2$ by using two three-manifolds $M_1$ and $M_2$ whose boundaries are tori, where we identify $\partial{M_1}$ with $-\partial{M_2}$.
We give the same basis for $\pi_1(\partial{M_1})$ and $\pi_1(-\partial{M_2})$ and let $\rho_i$ be the restriction to $\pi_1(M_i)$ of a representation $\rho$ of $\pi_1(W)$ ($i=1,2$).
Then we have
\begin{equation*}
  \cs_{W}\bigl([\rho]\bigr)
  =
  z_1z_2,
\end{equation*}
where $\cs_{M_i}\left(\left[\rho_i\right]\right)=[\alpha,\beta,z_i]$ with respect to the common basis \cite[Theorem~2.2]{Kirk/Klassen:COMMP93}.
\par
For a three-manifold $M$ with hyperbolic metric one can define the ($SO(3)$) Chern--Simons invariant $\CS(M)$ as follows \cite{Chern/Simons:ANNMA21974}.
Let $A$ be the Levi-Civita connection (an $so(3)$-valued $1$-form) defined by the hyperbolic metric.
Then put
\begin{equation*}
  \cs(M)
  :=
  \frac{1}{8\pi^2}
  \int_{M}\tr
  \left(A\wedge{dA}+\frac{2}{3}A\wedge A\wedge A\right)
  \in\R
  \pmod{\Z}.
\end{equation*}
In this paper we use another normalization $\CS(M):=-2\pi^2\cs(M)$ so that $\Vol(M)+\sqrt{-1}\CS(M)$ can be regarded as a complexification of the volume $\Vol(M)$.
Note that if $M$ has a cusp (that is, if $M$ is homeomorphic to the interior of a manifold with torus boundary), $\cs$ is defined modulo $1/2$ and so $\CS$ is defined modulo $\pi^2$.
It is known that the imaginary part of the $SL(2;\C)$ Chern--Simons function of the holonomy representation is the volume with respect to the hyperbolic metric and the real part is the $SO(3)$ Chern--Simons invariant up to multiplications of some constants \cite[Lemma~3.1]{Yoshida:INVEM85} (see also \cite[p.~554]{Kirk/Klassen:COMMP93}).
\par
When $M$ has a cusp and possesses a complete hyperbolic structure of finite volume, we can deform the structure to incomplete ones by using a complex parameter $u$ around $0$, where $u=0$ corresponds to the complete structure.
Then T.~Yoshida \cite[Theorem~2]{Yoshida:INVEM85} (see also \cite[Conjecture]{Neumann/Zagier:TOPOL85}) proved that there exists a complex analytic function $f(u)$ around $0$ such that if the corresponding (incomplete) hyperbolic manifold can be completed to a closed manifold $M_u$ by adding a geodesic loop $\gamma$, then
\begin{multline}\label{eq:Vol_CS_f}
  \Vol(M_u)+\sqrt{-1}\CS(M_u)
  -\left\{\Vol(M)+\sqrt{-1}\CS(M)\right\}
  \\
  \equiv
  \frac{f(u)}{\sqrt{-1}}
  -\frac{\pi}{2}
  \left\{\length(\gamma)+\sqrt{-1}\torsion(\gamma)\right\}
  \pmod{\pi^2\sqrt{-1}\Z},
\end{multline}
where $\length(\gamma)$ and $\torsion(\gamma)$ are the length and the torsion of $\gamma$ respectively.
(The torsion measures how much a normal vector is twisted when it travels around $\gamma$.)
\par
We can interpret \eqref{eq:Vol_CS_f} in terms of Kirk and Klassen's Chern--Simons function $\cs_M$.
\par
For simplicity we assume that $M:=S^3\setminus{K}$ is the complement of a hyperbolic knot $K$, that is, we assume that $M$ has a complete hyperbolic structure with finite volume.
Then $u$ defines an (incomplete) hyperbolic structure, and so gives the holonomy representation $\rho$ of $\pi_1(M)$ at $PSL(2;\C)$.
(Recall that the orientation preserving isometry group of the hyperbolic space $\mathbf{H}^3$ is $PSL(2;\C)$.
Since the universal cover of $M$ is $\mathbf{H}^3$, the lift of an element in $\pi_1(M)$ defines an isometric translation in $\mathbf{H}^3$, giving an element in $PSL(2;\C)$.)
We choose a lift of $\rho$ to $SL(2;\C)$ \cite{Culler/ADVAM41984} and denote it also by $\rho$.
\par
We can regard $\exp(u)$ as the ratio of the eigenvalues of the image of the meridian $\mu$ by the representation $\rho$.
We can also define $v(u)$ so that $\exp\bigl(v(u)\bigr)$ is the ratio of the eigenvalues of the image of the longitude $\lambda$.
Then from \cite[pp.~553--556]{Kirk/Klassen:COMMP93} $\cs_{K}:=\cs_{M}$ can be expressed in terms of $f(u)$ as follows.
\begin{equation}\label{eq:cs_f}
  \cs_{K}\bigl([\rho]\bigr)
  =
  \left[
    \frac{u}{4\pi\sqrt{-1}},\frac{v(u)}{4\pi\sqrt{-1}};
    \exp\left(\frac{\sqrt{-1}}{2\pi}f(u)\right)
  \right].
\end{equation}
Note that here we use $(\mu,\lambda)$ for the basis of $\pi_1(\partial{M})$.
\begin{rem}
Note that here we are using the $SL(2;\C)$ theory; not $PSL(2;\C)$.
So we have to divide the term appearing in $\exp$ of \eqref{eq:cs_f} by $-4$ \cite[pp.~553--556]{Kirk/Klassen:COMMP93} .
See \cite[p.~543]{Kirk/Klassen:COMMP93}.
\par
Note also that we are using the normalization of W.~Neumann and D.~Zagier~\cite{Neumann/Zagier:TOPOL85} for $f(u)$.
So Yoshida's (and Kirk and Klassen's) $f(u)$ is our $f(u)\times\dfrac{2}{\pi\sqrt{-1}}$.
\end{rem}
\par
If the incomplete metric is completed by adding a geodesic loop $\gamma$, then the resulting manifold $K_u$ is obtained from $S^3$ by $(p,q)$-Dehn surgery along $K$ for some coprime integers $p$ and $q$ satisfying $p\,u+q\,v(u)=2\pi\sqrt{-1}$.
Then from \cite[Lemma~4.2]{Neumann/Zagier:TOPOL85} we have
\begin{equation*}
  \length(\gamma)+\sqrt{-1}\torsion(\gamma)
  =
  -r\,u-s\,v(u)
  \pmod{2\pi\sqrt{-1}\Z}.
\end{equation*}
where $r$ and $s$ are integers satisfying $ps-qr=1$.
\par
Therefore we can express the right hand side of \eqref{eq:Vol_CS_f} in terms of $u$ and $v(u)$.
\begin{multline}\label{eq:Vol_CS_f_uv}
  \Vol(K_u)+\sqrt{-1}\CS(K_u)
  -\left\{\Vol(K)+\sqrt{-1}\CS(K)\right\}
  \\
  \equiv
  \frac{f(u)}{\sqrt{-1}}
  +\frac{\pi}{2}\bigl(r\,u+s\,v(u)\bigr)
  \pmod{\pi^2\sqrt{-1}\Z},
\end{multline}
where $\Vol(K):=\Vol(S^3\setminus{K})$ and $\CS(K):=\CS(S^3\setminus{K})$.
Note that this formula can also be obtained from \eqref{eq:cs_f} \cite{Kirk/Klassen:COMMP93}.
\par
We also have from \cite[(34), p.~323]{Neumann/Zagier:TOPOL85}
\begin{equation*}
  \length(u)
  =
  -\frac{1}{2\pi}\Im(u\overline{v(u)}).
\end{equation*}
\par
Therefore from \eqref{eq:Vol_CS_f} and \eqref{eq:Vol_CS_f_uv} we have
\begin{align}
  \Vol(K_u)-\Vol(K)
  &=
  \Im{f(u)}
  +\frac{1}{4}\Im(u\overline{v(u)}),
  \label{eq:Vol_CS_f_uv1}
  \\
  \CS(K_u)-\CS(K)
  &\equiv
  -\Re{f(u)}
  +\frac{\pi}{2}\Im(r\,u+s\,v(u))
  \pmod{\pi^2\Z}.
  \label{eq:Vol_CS_f_uv2}
\end{align}
Note that unfortunately we cannot express the Chern--Simons invariant only in terms of $u$.
%%%%%%%%%%%%%%%%%%%%%%%%%%%%%%%%%%%%%%%%%%%%%%%%%%%%%%%%%%%%%%%%%%%%%
\section{Volume Conjecture and its generalizations}
In \cite{Kashaev:MODPLA95} Kashaev introduced link invariants parametrized by an integer $N$ greater than one, by using the quantum dilogarithm.
Moreover he observed in \cite{Kashaev:LETMP97} that for the hyperbolic knots $4_1$, $5_2$, and $6_1$ his invariants grow exponentially and the growth rates give the hyperbolic volumes of the knot complements.
He also conjectured this would also hold for any hyperbolic knot in $S^3$.
\par
In \cite{Murakami/Murakami:ACTAM12001} J.~Murakami and the author showed that Kashaev's invariant coincides with the $N$-colored Jones polynomial evaluated at the $N$-th root of unity and generalized his conjecture to general knots as follows.
\begin{conj}[Volume Conjecture]\label{conj:VC}
For any knot $K$ in $S^3$ we would have
\begin{equation}\label{eq:VC}
  2\pi\lim_{N\to\infty}\frac{\log\left|J_N(K;\exp(2\pi\sqrt{-1}/N))\right|}{N}
  =
  \Vol(K),
\end{equation}
where $\Vol(K)$ is the simplicial volume {\rm(}or Gromov norm{\rm)} of the knot complement $S^3\setminus{K}$ {\rm(}see \cite{Gromov:INSHE82} and \cite[Chapter~6]{Thurston:GT3M}{\rm)}.
Note that $\Vol(K)$ is normalized so that it coincides with the hyperbolic volume if $K$ is a hyperbolic knot.
\end{conj}
This conjecture is so far proved for the following knots and links:
\begin{itemize}
\item
torus knots by Kashaev and O.~Tirkkonen \cite{Kashaev/Tirkkonen:ZAPNS2000}.
Note that their simplicial volumes are zero since their complements contain no hyperbolic pieces.
So in fact they proved that the left hand side of \eqref{eq:VC} converges to $0$.
See \cite{Dubois/Kashaev:MATHA2007} for more precise asymptotic behaviors.
See also \cite{Hikami:EXPMA2003,Hikami:COMMP2004,Hikami:RAMJO2006,Hikami/Kirillov:PHYLB2003} for related topics.
\item
the figure-eight knot $4_1$ by T.~Ekholm (see for example \cite[Section~3]{Murakami:Novosibirsk}),
\item
the hyperbolic knots $5_2$, $6_1$, and $6_2$ by Y.~Yokota,
\item
Whitehead doubles of torus knots of type $(2,b)$ by H.~Zheng \cite{Zheng:CHIAM22007},
\item
twisted Whitehead links by Zheng \cite{Zheng:CHIAM22007},
\item
the Borromean rings by S.~Garoufalidis and L{\^e} \cite[Theorem~12]{Garoufalidis/Le:2005},
\item
Whitehead chains, which generalizes both the Borromean rings and twisted Whitehead links, by R.~van der Veen \cite{van_der_Veen:2006}.
\end{itemize}
\par
Removing the absolute value sign of the left hand side of \eqref{eq:VC} we expect the Chern--Simons invariant as its imaginary part.
\begin{conj}[Complexification of the Volume Conjecture, \cite{Murakami/Murakami/Okamoto/Takata/Yokota:EXPMA02}]
For any knot $K$ in $S^3$ we would have
\begin{equation*}
  2\pi\lim_{N\to\infty}\frac{\log {J_N(K;\exp(2\pi\sqrt{-1}/N))}}{N}
  =
  \Vol(K)+\sqrt{-1}\CS(K),
\end{equation*}
where $\CS(K)$ is the Chern--Simons invariant of the knot complement \cite{Meyerhoff:density} if $K$ is a hyperbolic knot.
Note that we may regard the imaginary part of the left hand side as a definition of a {\it topological} Chern--Simons invariant.
\end{conj}
In \cite{Murakami/Murakami/Okamoto/Takata/Yokota:EXPMA02}, J.~Murakami, M.~Okamoto, T.~Takata, Y.~Yokota and the author used computer to confirm the Complexification of the Volume Conjecture for hyperbolic knots $6_3$, $8_9$, $8_{20}$ and for the Whitehead link up to several digits.
(So we do not have rigorous proofs.)
\par
What happens if we replace $2\pi\sqrt{-1}$ with another complex number?
Here is our conjecture which generalizes the complexification above.
%%%%%%%%%%%%%%%
\begin{conj}\label{conj:new}
Let $K$ be a knot in $S^3$.
Then there exists an open set $\mathcal{O}$ of $\C$ such that for any $u\in\mathcal{O}$ the series
$
\left\{
  \log
  \left(
    J_N
    \left(
      K;\exp\bigl((u+2\pi\sqrt{-1})/N\bigr)
    \right)
  \right)/N
\right\}_{N=2,3,\dots}
$
converges.
Moreover if we put
\begin{align}\label{eq:def_H_v_f}
  H(u)
  &:=
  (u+2\pi\sqrt{-1})
  \lim_{N\to\infty}
  \frac{\log{J_N\left(K;\exp\bigl((u+2\pi\sqrt{-1})/N\bigr)\right)}}{N},
  \\
  v(u)
  &:=
  2\frac{d\,H(u)}{d\,u}-2\pi\sqrt{-1},
  \\
  \intertext{and}
  f(u)
  &:=
  H(u)-H(0)-\pi\sqrt{-1}\,u-\frac{1}{4}uv(u),
\end{align}
then $f(u)$ becomes the $f$ function appearing in \eqref{eq:Vol_CS_f} and \eqref{eq:cs_f} in Subsection~\ref{subsec:V_CS} modulo $\pi^2$ and $H(0)$.
\par
Note that $H(0)$ may not be defined.
In that case $f(u)$ is defined modulo $H(0)$.
\end{conj}
Assuming Conjecture~\ref{conj:new} above we have from \eqref{eq:Vol_CS_f_uv1}
\begin{equation*}
\begin{split}
  \Vol(K_u)-\Vol(K)
  &=
  \Im{H(u)}-\Im{H(0)}-\pi\Re{u}-\frac{1}{4}\Im{uv(u)}
  +\frac{1}{4}\Im{u\overline{v(u)}}
  \\
  &=
  \Im{H(u)}-\Im{H(0)}-\pi\Re{u}-\frac{1}{2}\Re(u)\Im\bigl(v(u)\bigr).
\end{split}
\end{equation*}
From the Volume Conjecture \eqref{eq:VC} this is {\it almost} the same as the Parametrized Volume Conjecture \cite[Conjecture~2.1]{Murakami:ADVAM22007}, but note that we do not assume that the open set $\mathcal{O}$ contains $0$.
\par
In the following two sections we will show supporting evidence for Conjecture~\ref{conj:new} above, giving the figure-eight knot and torus knots as examples.
%%%%%%%%%%%%%%%%%%%%%%%%%%%%%%%%%%%%%%%%%%%%%%%%%%%%%%%%%%%%%%%%%
\section{Example~1 -- figure-eight knot}
In this section we describe how the colored Jones polynomials of the figure-eight knot $E$ are related to representations of the fundamental group at $SL(2;\C)$ and the corresponding volumes and Chern--Simons invariants.
%%%%%%%%%%%%%%%%%%%%%%%%%%%%%
\subsection{Representations of the fundamental group}\label{subsec:rep_fig8}
\begin{figure}[h]
\includegraphics[scale=0.5]{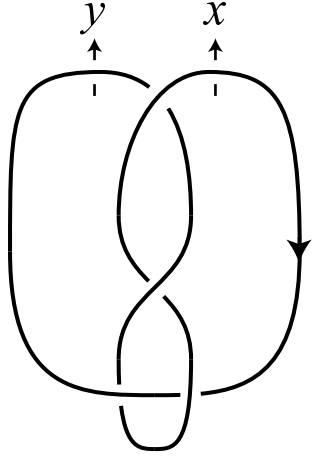}
\caption{The two generators of $\pi_1(S^3\setminus{E})$.}
\label{fig:fig8_pi1}
\end{figure}

\begin{figure}[h]
\includegraphics[scale=0.5]{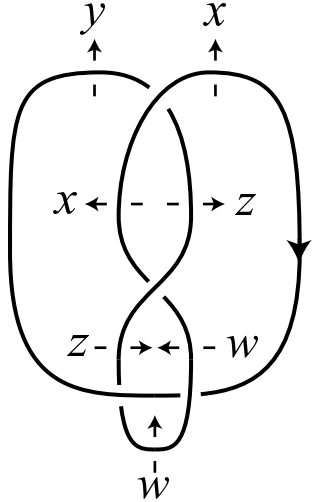}
\caption{The other elements of $\pi_1(S^3\setminus{E})$.}
\label{fig:fig8_pi1_relation}
\end{figure}

Here we follow R.~Riley \cite{Riley:QUAJM31984} to describe non-abelian representations of $\pi_1(S^3\setminus{E})$ at $SL(2;\C)$, where $E$ is the figure-eight knot.
Let $x$ and $y$ be the generators of $\pi_1(S^3\setminus{E})$ indicated in Figure~\ref{fig:fig8_pi1}, where the basepoint of the fundamental group is above the paper.
%%%%%%%%%%%%
%%%%%%%%%%%%
Then the other elements $z$ and $w$ indicated in Figure~\ref{fig:fig8_pi1_relation} can be expressed by $x$ and $y$ as follows.
%%%%%%%%%%%%%%
%%%%%%%%%%%%%
\begin{align}
  z&=xyx^{-1}
  \label{eq:fig8_z}
  \\
  &\quad\text{(from the top-most crossing)},
  \notag
  \\
  w&=z^{-1}xz
  \notag
  \\
  &\quad\text{(from the second top crossing)}.
  \notag
\end{align}
Therefore we have
\begin{equation}\label{eq:fig8_w}
  w=\left(xy^{-1}x^{-1}\right)x\left(xyx^{-1}\right)
   =xy^{-1}xyx^{-1}
\end{equation}
Now from the bottom-left and the bottom-right crossings we have the following relations.
\begin{align*}
  yw^{-1}y^{-1}z&=1,
  \\
  yw^{-1}x^{-1}w&=1.
\end{align*}
Using \eqref{eq:fig8_z} and \eqref{eq:fig8_w}, these relations become a single relation
\begin{equation}\label{eq:fig8_relation}
  y
  xy^{-1}x^{-1}yx^{-1}
  y^{-1}
  xyx^{-1}
  =1.
\end{equation}
Therefore $\pi_1(S^3\setminus{E})$ has the two generators $x$ and $y$ with the single relation \eqref{eq:fig8_relation}, that is:
\begin{equation}\label{eq:fig8_presentation}
  \pi_1(S^3\setminus{E})
  =
  \langle
    x,y
    \mid
    \omega x=y\omega
  \rangle
\end{equation}
with $\omega:=xy^{-1}x^{-1}y$.
\par
Let $\rho$ be a non-abelian representation of $\pi_1(S^3\setminus{E})$ at $SL(2;\C)$.
From \cite{Riley:QUAJM31984} we can assume up to conjugate that the images of $\rho(x)$ and $\rho(y)$ are given as follows.
\begin{align*}
  \rho(x)
  &:=
  \begin{pmatrix}
    m^{1/2} & 1 \\
    0       & m^{-1/2}
  \end{pmatrix},
  \\
  \rho(y)
  &:=
  \begin{pmatrix}
    m^{1/2} & 0 \\
    -d      & m^{-1/2}
  \end{pmatrix}.
\end{align*}
Then from the presentation \eqref{eq:fig8_presentation}, we have $\rho(\omega x)=\rho(y\omega)$.
Since $\rho(\omega)$ is equal to
\begin{equation*}
  \begin{pmatrix}
    (d+1)^2-dm & -m^{-1/2}(d-m+1)
    \\
    dm^{-1/2}(d-m+1) & m^{-1}(m-d)
  \end{pmatrix},
\end{equation*}
we have
\begin{equation*}
  \rho(\omega x)
  =
  \begin{pmatrix}
    m^{1/2}((d+1)^2-dm) & -m^{-1}(dm^2-(d^2+2d+2)m+d+1) \\
    -d(m-d-1) & m^{-3/2}(-dm^2+(d^2+d+1)m-d)
  \end{pmatrix},
\end{equation*}
and
\begin{multline*}
  \rho(y\omega)
  \\
  =
  \begin{pmatrix}
    m^{1/2}((d+1)^2-dm) & m-d-1 \\
    dm^{-1}(dm^2-(d^2+2d+2)m+d+1) & m^{-3/2}(-dm^2+(d^2+d+1)m-d)
  \end{pmatrix}.
\end{multline*}
So $d$ and $m$ satisfy the equation
\begin{equation}\label{eq:d_m_fig8}
  d^2+d(3-m-m^{-1})+3-m-m^{-1}=0,
\end{equation}
and $d$ becomes a function of $m$:
\begin{equation*}
  d
  =
  \frac{1}{2}
  \left(m+m^{-1}-3\pm\sqrt{(m+m^{-1}+1)(m+m^{-1}-3)}\right).
\end{equation*}
Therefore if we put $\xi:=\tr\bigl(\rho(x)\bigr)$ and $\eta:=\tr\bigl(\rho(xy)\bigr)$ as in Subsection~\ref{subsec:character_variety}, we have
\begin{align*}
  \xi
  &=
  m^{1/2}+m^{-1/2}
  \\
  \eta
  &=
  \tr
  \begin{pmatrix}
    m-d & m^{-1/2} \\
    -dm^{-1/2} & m^{-1}
  \end{pmatrix}
  =
  m+m^{-1}-d.
\end{align*}
Moreover since from \eqref{eq:trace} we calculate
\begin{align*}
  \tr\bigl(\omega\bigr)&=\eta^2-\xi^2\eta+2\xi^2-2
  \\
  \tr\bigl(\omega^{(1)}\bigr)&=\eta
\end{align*}
we have
\begin{equation*}
  F(\xi,\eta)=\eta^2-\eta+2\xi^2-\xi^2\eta-1,
\end{equation*}
which also gives \eqref{eq:d_m_fig8}.
(See Theorem~\ref{thm:Le}.)
\par
Now the longitude $\lambda$ is (read off from the top right)
\begin{equation*}
  \lambda
  =
  wx^{-1}yz^{-1}
  =
  xy^{-1}xyx^{-1} x^{-1} y xy^{-1}x^{-1}
  =
  xy^{-1}xyx^{-2}yxy^{-1}x^{-1}
\end{equation*}
From a direct calculation, we have
\begin{equation*}
  \rho(\lambda)
  =
  \begin{pmatrix}
    \ell(m)^{\pm1}
    &
    (m^{1/2}+m^{-1/2})\sqrt{(m+m^{-1}+1)(m+m^{-1}-3)}
    \\
    0
    &
    \ell(m)^{\mp1}
  \end{pmatrix},
\end{equation*}
where
\begin{multline*}
  \ell(m)
  \\
  :=
  \frac{\left(m^2-m-2-m^{-1}+m^{-2}\right)}{2}
  +\frac{(m-m^{-1})}{2}
  \sqrt{(m+m^{-1}+1)(m+m^{-1}-3)}.
\end{multline*}
Note that $\ell(m)$ is a solution to the following equation.
\begin{equation}\label{eq:A_polynomial}
  \ell+(m^2-m-2-m^{-1}+m^{-2})+\ell^{-1}=0,
\end{equation}
which coincides with the $A$-polynomial of the figure-eight knot \cite{Cooper/Culler/Gillet/Long/Shalen:INVEM1994}, replacing $\ell$ and $m$ with $\mathfrak{l}$ and $\mathfrak{m}^2$, respectively.
\par
Let us denote the representations $\rho_{m+}$ and $\rho_{m-}$.
%%%%%%%%%%%%%%%%%%%%%%%%%%%%%%%%%%%%%%%%%%%%%%%%%%%%%%
\subsection{Asymptotic behavior of the colored Jones polynomial}
Now we consider the colored Jones polynomial and its asymptotic behavior.
\par
As described above, Habiro and L{\^e} independently showed that for the figure-eight knot $E$, $J_N(E;q)$ can be expressed in a single summation as follows:
\begin{equation}\label{eq:J_N(E;q)}
  J_N(E;q)
  =
  \sum_{k=0}^{N-1}
  \prod_{j=1}^{k}
  \left(q^{(N-j)/2}-q^{-(N-j)/2}\right)
  \left(q^{(N+j)/2}-q^{-(N+j)/2}\right).
\end{equation}
\par
We put $q:=\exp(\theta/N)$ and consider the asymptotic behavior of $J_N(E;\exp(\theta/N))$ for $N\to\infty$ fixing a complex parameter $\theta$.
\par
To state a formula describing the asymptotic behavior of $J_N(E;\exp(\theta/N))$ we prepare some functions.
Put $\varphi(\theta):=\arccosh(\cosh(\theta)-1/2)$, where we choose the branch of $\arccosh$ so that
\begin{equation*}
  \arccosh(x)
  =
  \log\left(x-\sqrt{-1}\sqrt{1-x^2}\right)+2\pi\sqrt{-1}.
\end{equation*}
We also choose the branch cut of $\log$ as $(-\infty,0)$.
Note that we use $\sqrt{-1}\sqrt{1-x^2}$ instead of $\sqrt{x^2-1}$ to avoid the cut branch of the square root function, since later we will assume that $x$ is near to $1/2$.
Let $\Li_2$ be the dilogarithm function:
\begin{equation*}
  \Li_2(z)
  :=
  -\int_{0}^{z}\frac{\log(1-w)}{w}\,dw
  =
  \sum_{k=1}^{\infty}\frac{z^k}{k^2}
\end{equation*}
with the branch cut $(1,\infty)$.
\par
Now we can prove the following four theorems about the asymptotic behavior of the colored Jones polynomials of the figure-eight knot $E$.
\par\medskip\noindent
(i). When $\theta$ is close to $2\pi\sqrt{-1}$.
%%%%%%%%%%
\begin{thm}[\cite{Murakami/Yokota:JREIA2007}]\label{thm:MY}
Let $\theta$ be a complex number near to $2\pi\sqrt{-1}$.
We also assume that $\theta$ is not purely imaginary except for $\theta=2\pi\sqrt{-1}$.
Then we have
\begin{equation*}
  \theta
  \lim_{N\to\infty}
  \frac{\log{J_N\bigl(E;\exp(\theta/N)\bigr)}}{N}
  =
  H(\theta),
\end{equation*}
where we put
\begin{equation*}
  H(\theta)
  :=
  \Li_2\left(e^{-\varphi(\theta)-\theta}\right)
  -
  \Li_2\left(e^{\varphi(\theta)-\theta}\right)
  +(\theta-2\pi\sqrt{-1})\varphi(\theta).
\end{equation*}
\end{thm}
%%%%%%%%%%%%%%%%%%%%%%%%%%
\par\medskip\noindent
(ii). When $\theta$ is real and $|\theta|\ge\arccosh(3/2)$.
\par
\begin{thm}[{\cite[Theorem~8.1]{Murakami:KYUMJ2004} (See also \cite[Lemma~6.7]{Murakami:Novosibirsk})}]
Let $\theta$ be a real number with $|\theta|\ge\arccosh(3/2)$.
Then we have
\begin{equation*}
  \theta
  \lim_{N\to\infty}
  \frac{\log{J_N\bigl(E;\exp(\theta/N)\bigr)}}{N}
  =
  \tilde{H}(\theta),
\end{equation*}
where we put
\begin{equation*}
  \tilde{H}(\theta)
  :=
  \Li_2\left(e^{-\tilde{\varphi}(\theta)-\theta}\right)
  -
  \Li_2\left(e^{\tilde{\varphi}(\theta)-\theta}\right)
  +
  \theta\tilde{\varphi}(\theta)
\end{equation*}
and $\tilde{\varphi}(\theta)$ is defined as $\varphi(\theta)$ by using the usual branch of $\arccosh$ so that $\arccosh(x)>0$ if $|x|>1$.
\par
Note that $\tilde{H}(\theta)=0$ if $\theta=\pm\arccosh(3/2)$.
\end{thm}
\par\medskip\noindent
(iii) When $\theta$ is close to $0$.
\begin{thm}[\cite{Murakami:JPJGT2007},\cite{Garoufalidis/Le:aMMR}]\label{thm:MGL}
Let $\theta$ be a complex number with $|2\cosh(\theta)-2|<1$ and $\Im(\theta)<\pi/3$, then the series $\{J_N\bigl(E;\exp(\theta/N)\bigr)\}_{N=2,3,\dots}$ converges and
\begin{equation*}
  \lim_{N\to\infty}J_N\bigl(E;\exp(\theta/N)\bigr)
  =
  \frac{1}{\Delta\bigl(E;\exp(\theta)\bigr)},
\end{equation*}
where $\Delta(E;t)$ is the Alexander polynomial of the figure-eight knot $E$.
\par
Here we normalize the Alexander polynomial for a knot $K$ so that $\Delta(K;1)=1$ and that $\Delta(K;t^{-1})=\Delta(K;t)$.
\par
Note that in this case
\begin{equation*}
  \lim_{N\to\infty}
  \frac{\log{J_N\bigl(E;\exp(\theta/N)\bigr)}}{N}
  =0.
\end{equation*}
\end{thm}
\par\medskip\noindent
(iv). When $\theta=\pm\arccosh(3/2)$.
\par
From (ii) we know that when $\theta=\pm\arccosh(3/2)$ then the colored Jones polynomial does not grow exponentially.
Since  $\exp\bigl(\pm\arccosh(3/2)\bigr)=(3\pm\sqrt{5})/2$ is a zero of $\Delta(E;t)$, we cannot expect that the colored Jones polynomial converges.
In fact in this case we can prove that it grows polynomially.
\begin{thm}[{\cite[Theorem~1.1]{Hikami/Murakami:2007}}]
We have
\begin{equation*}
  J_N(E;\exp(\pm\arccosh(3/2)/N))
  \underset{N\to\infty}{\sim}
  \frac{\Gamma(1/3)}{\bigl(3\arccosh(3/2)\bigr)^{2/3}}
  N^{2/3},
\end{equation*}
where $\Gamma$ is the gamma function:
\begin{equation*}
  \Gamma(z)
  :=
  \int_{0}^{\infty}t^{z-1}e^{-t}\,dt.
\end{equation*}
\end{thm}
%%%%%%%%%%%%%%%%%%%%%%%%%%%%%%%%%%%%%%
\subsection{Relation of the function $H$ to a representation}
\label{subsec:H_rep_fig8}
Now we want to relate the function $H$ to a representation.
It is convenient to put $u:=\theta-2\pi\sqrt{-1}$ so that the complete hyperbolic structure corresponds to $u=0$.
\par\medskip\noindent
(i). When $\theta$ is close to $2\pi\sqrt{-1}$, that is, when $u$ is close to $0$.
\par
We put
\begin{equation}\label{eq:Hu}
  H(u)
  :=
  \Li_2\left(e^{-\varphi(u)-u}\right)
  -\Li_2\left(e^{\varphi(u)-u}\right)
  +u\varphi(u).
\end{equation}
Let us calculate the derivative $d\,H(u)/d\,u$.
\begin{comment}
Since
\begin{equation*}
  \frac{d\,\arccosh(x)}{d\,x}
  =
  \frac{1+\frac{\sqrt{-1}x}{\sqrt{1-x^2}}}{x-\sqrt{-1}\sqrt{1-x^2}}
  =
  \frac{\sqrt{-1}}{\sqrt{1-x^2}},
\end{equation*}
we have
\begin{equation*}
  \frac{d\,\varphi(u)}{d\,u}
  =
  \frac{4\sqrt{-1}\sinh(u)}{\sqrt{3+4\cosh(u)-4\cosh^2(u)}}.
\end{equation*}
\end{comment}
We have
\begin{equation*}
\begin{split}
  \frac{d\,H(u)}{d\,u}
  \\
  =&
  -
  \frac{\log\left(1-e^{-\varphi(u)-u}\right)}{e^{-\varphi(u)-u}}
  \times
  e^{-\varphi(u)-u}
  \times
  \left(-\frac{d\,\varphi(u)}{d\,u}-1\right)
  \\
  &+
  \frac{\log\left(1-e^{\varphi(u)-u}\right)}{e^{\varphi(u)-u}}
  \times
  e^{\varphi(u)-u}
  \times
  \left(\frac{d\,\varphi(u)}{d\,u}-1\right)
  +
  \varphi(u)
  +
  u\times\frac{d\,\varphi(u)}{d\,u}
  \\
  =&
  \frac{d\,\varphi(u)}{d\,u}
  \left\{
    \log\left(1-e^{-\varphi(u)-u}\right)
    +
    \log\left(1-e^{\varphi(u)-u}\right)
    +
    u
  \right\}
  \\
  &+
  \log\left(1-e^{-\varphi(u)-u}\right)
  -
  \log\left(1-e^{\varphi(u)-u}\right)
  +
  \varphi(u).
\end{split}
\end{equation*}
Since $\varphi(u)$ satisfies $\cosh\bigl(\varphi(u)\bigr)=\cosh(u)-1/2$, we have
\begin{equation*}
  e^{\varphi(u)}+e^{\varphi(-u)}=e^u+e^{-u}-1
\end{equation*}
and so
\begin{equation*}
  \left(1-e^{-\varphi(u)-u}\right)
  \left(1-e^{\varphi(u)-u}\right)
  =
  e^{-u}\left(e^u-e^{-\varphi(u)}-e^{\varphi(u)}+e^{-u}\right)
  =
  e^{-u}.
\end{equation*}
Therefore we have
\begin{equation}\label{eq:derivative_H}
\begin{split}
  \frac{d\,H(u)}{d\,u}
  &=
  2\log\left(1-e^{-\varphi(u)-u}\right)
  +
  \varphi(u)
  +
  u.
\end{split}
\end{equation}
\par
If we put
\begin{equation}\label{eq:v}
\begin{split}
  v(u)
  &:=
  2\frac{d\,H(u)}{d\,u}-2\pi\sqrt{-1}
  \\
  &=
  4\log\left(1-e^{-\varphi(u)-u}\right)
  +
  2\varphi(u)
  +
  2u
  -2\pi\sqrt{-1},
\end{split}
\end{equation}
we have
\begin{equation*}
  \exp\left(\frac{v(u)}{2}\right)
  =
  -
  \exp
  \left(
    \log
    \left(
      e^{\varphi(u)+u}-2+e^{-\varphi(u)-u}
    \right)
  \right)
  =
  -e^{\varphi(u)+u}+2-e^{-\varphi(u)-u}.
\end{equation*}
We also put $m:=\exp(u)$.
Then since
\begin{equation*}
\begin{split}
  e^{\pm\varphi(u)}
  &=
  \cosh(u)-\frac{1}{2}
  \mp
  \sqrt{-1}\sqrt{1-\left(\cosh(u)-\frac{1}{2}\right)^2}
  \\
  &=
  \frac{1}{2}\left(m+m^{-1}-1\right)
  \mp
  \frac{\sqrt{-1}}{2}
  \sqrt{(m+1+m^{-1})(3-m-m^{-1})},
\end{split}
\end{equation*}
we have
\begin{equation*}
\begin{split}
  \exp\left(\frac{v(u)}{2}\right)
  =&
  \frac{1}{2}
  \left(
    -m^2+m+2+m^{-1}-m^{-2}
  \right)
  \\
  &-
  \frac{\sqrt{-1}\left(m-m^{-1}\right)}{2}
  \sqrt{(m+1+m^{-1})(3-m-m^{-1})}
  \\
  =&
  -\ell(m).
\end{split}
\end{equation*}
\par
Therefore the representation $\rho_{m\pm}$ sends the longitude $\lambda$ to
\begin{equation*}
  \begin{pmatrix}
    -e^{\pm v(u)/2} & \ast \\
       0        & -e^{\mp v(u)/2}
  \end{pmatrix}.
\end{equation*}
\par\medskip\noindent
(ii). When $\theta$ is real and $|\theta|\ge\arccosh(3/2)$, that is, $\Im{u}=-2\pi\sqrt{-1}$ and $|\Re{u}|\ge\arccosh(3/2)$.
\par
In this case we put
\begin{equation*}
  \tilde{H}(u)
  :=
  \Li_2\left(e^{-\tilde{\varphi}(u)-u}\right)
  -
  \Li_2\left(e^{\tilde{\varphi}(u)-u}\right)
  +(u+2\pi\sqrt{-1})\tilde{\varphi}(u).
\end{equation*}
If we define $v(u)$ as in the previous case by using $\tilde{H}$ instead of $H$, then we see that $u$ defines the same representations as (i).
\par\medskip\noindent
(iii). When $\theta$ is close to $0$, that is $u$ is close to $-2\pi\sqrt{-1}$.
\par
In this case we may say that $\theta$ defines the abelian representation $\alpha_{\exp(\theta/2)}$ that sends the meridian element to the matrix
\begin{equation*}
  \begin{pmatrix}
    e^{\theta/2} & 0\\
      0     & e^{-\theta/2}
  \end{pmatrix}.
\end{equation*}
This is because the complex value $\Delta(E;\exp(\theta))$ is the determinant of the Fox matrix corresponding to the map $\pi_1(S^3\setminus{E})\to H_1(S^3;\Z)\cong\Z\to\C^{\ast}$, where the first map is the abelianization, and the second map sends $k$ to $\exp(k\theta)$.
See for example \cite[Chapter~11]{Lickorish:1997}.
\par\medskip\noindent
(iv). When $\theta=\pm\arccosh(3/2)$, that is, $u=\pm\arccosh(3/2)-2\pi\sqrt{-1}$.
\par
Note that when $u=\pm\arccosh(3/2)-2\pi\sqrt{-1}$, the corresponding representation is non-abelian from (ii), but `attached' to an abelian one from (iii).
\par
Let us study $\rho_{m\pm}$ defined in Subsection~\ref{subsec:rep_fig8} more carefully when $m=\exp\bigl(\pm\arccosh(3/2)\bigr)=(3\pm\sqrt{5})/2$.
\par
First note that $d=0$, and so $\rho_{m+}$ and $\rho_{m-}$ coincide.
Moreover the images of both $x$ and $y$ are upper triangle.
Therefore these representations are reducible.
\par
By the linear fractional transformation, the Lie group $SL(2;\C)$ acts on $\C$ by $\begin{pmatrix}p & q \\ r & s \end{pmatrix}\cdot z:=(pz+q)/(rz+s)$.
For an upper triangle matrix this becomes an affine transformation.
Therefore if $\theta=\pm\arccosh(3/2)$ the corresponding representations can be regarded as affine.
\par
It is known that each affine representation of a knot group corresponds to a zero of the Alexander polynomial of the knot (\cite{Burde:MATHA1967,deRham:ENSEM21967}).
See also \cite[Exercise~11.2]{Kauffman:Knots}.
\par
Note also that $\tr\bigl(\rho_{m\pm}\bigr)=\tr\bigl(\alpha_{m^{\pm1/2}}\bigr)$ in this case.
%%%%%%%%%%%%%%%%%%%%%%%%%%%%%%%%%%%%%%%%%%%%%%%%%%
\subsection{Volume and the Chern--Simons invariant}
\label{subsec:V_CS_fig8}
Now we know that $u$ defines a representation (with some ambiguity) of $\pi_1(S^3\setminus{E})$ at $SL(2;\C)$.
In this subsection we calculate the corresponding volume and Chern--Simons invariant.
\par\medskip\noindent
(i). When $u$ is close to $0$.
\par
In this case the parameter $u$ defines an (incomplete) hyperbolic structure on $S^3\setminus{E}$, and as described in Subsection~\ref{subsec:V_CS}, if $u$ and $v(u)$ satisfy $p\,u+q\,v(u)=2\pi\sqrt{-1}$ with coprime integers $p$ and $q$, then we can construct a closed three-manifold $E_u$.
In \cite{Murakami/Yokota:JREIA2007} we proved that Conjecture~\ref{conj:new} holds in this case.
Note that since $H(0)=\Vol(E)\sqrt{-1}$, which is the original volume conjecture \ref{conj:VC} for the figure-eight knot, $f(u)$ defined by $H(u)$ coincides with the $f$ function appearing in \eqref{eq:Vol_CS_f} and \eqref{eq:cs_f} without ambiguity.
%%%%%%%%%%%%%%%%%%%%%%%%%%%%%%%%%%
\begin{ex}\label{ex:real_u_volume}
As an example let us consider the case where $u$ is real, and study $H(u)$ and the corresponding representation.
\begin{rem}
In this example we do not mind whether the series
$\log\left(J_N\left(E;\exp\bigl((u+2\pi\sqrt{-1})/N\bigr)\right)\right)/N$ really converges or not.
\end{rem}
First we assume that $-\arccosh(3/2)<u<\arccosh(3/2)$.
Note that $\arccosh(3/2)=0.9624236501\ldots$.
\par
Since $|x-\sqrt{-1}\sqrt{1-x^2}|=1$ when $|x|\le1$ and $1/2\le\cosh(u)-1/2<1$ for real $u$ with $|u|\le\arccosh(3/2)$, we see that $\varphi(u)$ is purely imaginary.
Then from \eqref{eq:Hu} we have
\begin{equation*}
  H(u)
  =
  \Li_2(e^{-\varphi(u)-u})-\overline{\Li_2(e^{-\varphi(u)-u})}
  +u\varphi(u)
\end{equation*}
and so $H(u)$ is purely imaginary, where $\overline{z}$ is the complex conjugate of $z$.
We can also see that $v(u)$ is also purely imaginary from \eqref{eq:v}.
\par
Let us consider the corresponding volume function:
\begin{equation*}
  \Vol(E_u):=\Im{H(u)}-\pi\Re(u)-\frac{1}{2}\Re(u)\Im\bigl(v(u)\bigr).
\end{equation*}
\par
See Figures~\ref{fig:ImHu} and \ref{fig:vol} for graphs of $\Im{H(u)}$ and $\Vol(E_u)$ respectively.
%%%%%%%%%%%%%%%%%%
\begin{figure}[h]
\includegraphics[scale=1]{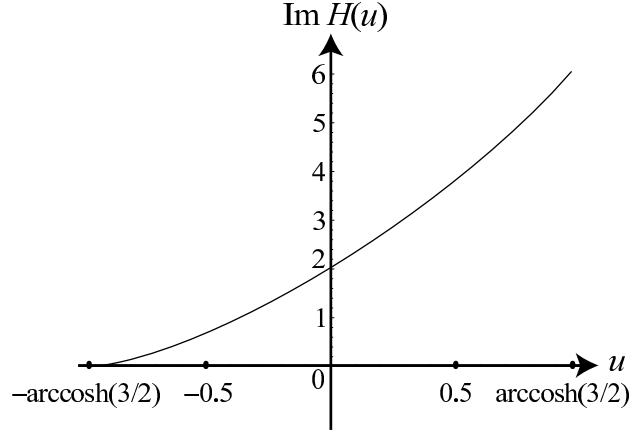}
\caption{Graph of $\Im{H(u)}$ for $-\arccosh(3/2)\le{u}\le\arccosh(3/2)$.}
\label{fig:ImHu}
\end{figure}
%%%%%%%%%%%%%%%%%%
\begin{figure}[h]
\includegraphics[scale=1]{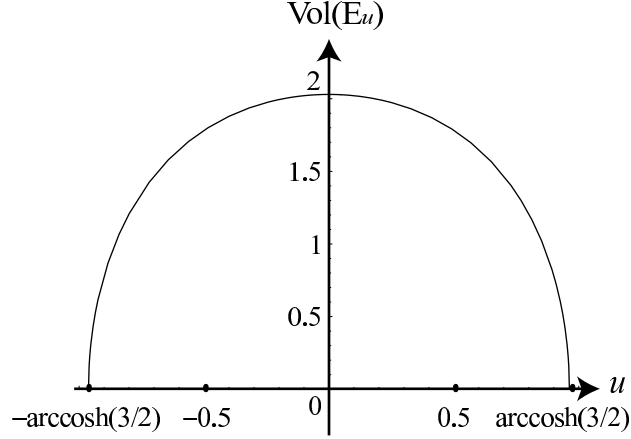}
\caption{Graph of the volume function $\Vol(E_u)$ for $-\arccosh(3/2)\le{u}\le\arccosh(3/2)$.}
\label{fig:vol}
\end{figure}
%%%%%%%%%%%%%%%%%%
Note that $\Vol(E_0)=\Im{H(0)}=2.0298832128\ldots$, which is the volume of $S^3\setminus{E}$ with the complete hyperbolic structure.
Note also that $\Im{H(\arccosh(3/2))}=2\pi\arccosh(3/2)=6.0470861377\ldots$ and $\Vol\left(E_{\pm\arccosh(3/2)}\right)=0$ since $\varphi(\arccosh(3/2))=2\pi\sqrt{-1}$.
\par
If we choose $u_q$ so that $v(u_q)=2\pi\sqrt{-1}/q$ for a positive integer $q$, then $E_{u_q}$ is the cone-manifold with cone-angle $2\pi/q$ whose underlying manifold is obtained from the $0$-surgery along the figure-eight knot (the singular set is the core of the surgery).
This is because it corresponds to the generalized Dehn surgery of coefficient $(0,q)$, since $0\times u_q+q\times v(u_q)=2\pi\sqrt{-1}$ \cite[Chapter~4]{Thurston:GT3M}.
Note that even though $q$ is not an integer, $\Vol(E_{u_q})$ is still the volume of the corresponding incomplete hyperbolic manifold.
\par
See \cite{Hilden/Lozano/Montesinos:JMASU1995} for geometric structures of the $E_{u_q}$.
They also observed that when $q=1$, that is, when $u=\pm\arccosh(3/2)$, the corresponding manifold is just the $0$-surgery and has a Sol-geometry.
They also calculated $\Vol(E_{u_q})$ with $q=2,3,\dots,81$ \cite[p.~559]{Hilden/Lozano/Montesinos:JMASU1995}.
%%%%%%%%%%%%%%%%%%%%%%%
\par
Next we consider the case where $u$ is real and $|u|>\arccosh(3/2)$.
In this case $\Im{v(u)}=2\pi\sqrt{-1}$ since $\Im{\varphi(u)}=2\pi\sqrt{-1}$.
\par
If we put $p=-\Re{v(u)}/u$, $q=1$, $s=0$ and $r=-1$, then $p\,u+q\,v(u)=2\pi\sqrt{-1}$ and $ps-qr=1$.
(For a while we do not mind whether $p$ is an integer or not.
See Remark~\ref{rem:Mathematica} below.)
Assuming that this would give a genuine manifold, we have from \eqref{eq:Vol_CS_f_uv} and \eqref{eq:def_H_v_f}
\begin{multline*}
  \Vol(E_u)+\sqrt{-1}\CS(M_u)-
  \left\{\Vol(E)+\sqrt{-1}\CS(E)\right\}
  \equiv
  \frac{H(u)}{\sqrt{-1}}
  -\frac{H(0)}{\sqrt{-1}}
  -\frac{3}{2}\pi{u}
  -\frac{uv(u)}{4\sqrt{-1}}
  \\
  \pmod{\pi^2\sqrt{-1}\Z}.
\end{multline*}
Since $H(0)=\sqrt{-1}\Vol(E)$ and $\CS(E)=0$ from the amphicheirality of $E$ (that is, the figure-eight knot is equivalent to its mirror image), we have
\begin{align*}
  \Vol(E_u)
  &=
  \Im{H(u)}-2\pi{u}
  =0
  \\
\intertext{and}
  \CS(E_u)
  &\equiv
  -\Re{H(u)}
  +\frac{1}{4}u\Re{v(u)}
  \pmod{\pi^2\Z},
\end{align*}
since $\Im{\varphi(u)}=2\pi$ if $|u|>\arccosh(3/2)$.
%%%%%%%%%%%%%%%%%%
\begin{figure}[h]
\includegraphics[scale=1]{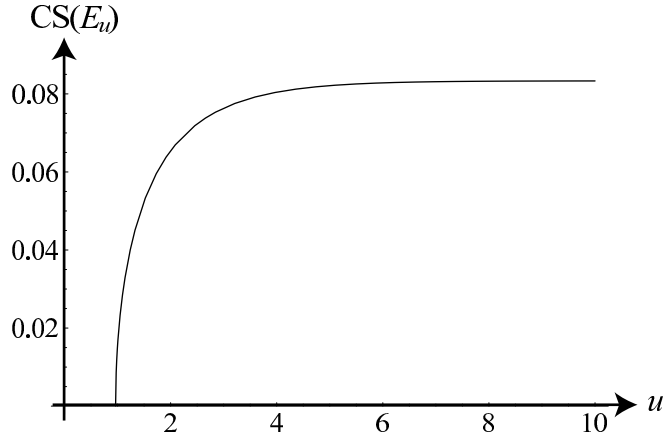}
\caption{Graph of the $\CS(M_u)/(2\pi^2)$ for $\le\arccosh(3/2)<u<10$.}
\label{fig:CS}
\end{figure}
%%%%%%%%%%%%%%%%%%
\end{ex}
%%%%%%%%%%%%%%%%%
\begin{rem}\label{rem:Mathematica}
Here we show some calculations by using Mathematica and SnapPea.
See Figure~\ref{fig:CS}.
\par
If $e^u$ satisfies the equation
\begin{equation*}
  x^6-4x^4-7x^3-4x^2+1=0
\end{equation*}
($u=0.9839865622\ldots$), then $u+v(u)=2\pi\sqrt{-1}$ and so this corresponds to the $(1,1)$-surgery along $E$.
Moreover we have $\dfrac{\CS(E_u)}{2\pi\sqrt{-1}}=0.01190476190\ldots=\dfrac{1}{84}$, which coincides with the calculation by SnapPea.
Note that $E_{u}$ in this case is the Brieskorn homology sphere $\Sigma(2,3,7)$, which is a Seifert fibered space over the sphere with three singular fibers.
(See for example \cite[Problem~1.77]{Kirby:problems}.)
\par
If
$u=
\log\left(
  \dfrac{1+\sqrt{5}}{2}+\sqrt{\dfrac{\sqrt{5}+1}{2}}
\right)
=1.061275062\ldots$,
then $2u+v(u)=2\pi\sqrt{-1}$ and so this corresponds to the $(2,1)$-surgery along $E$.
We also have $\dfrac{\CS(E_u)}{2\pi^2}=0.025=\dfrac{1}{40}$, which coincides with the calculation by SnapPea.
In this case $E_{u}$ is the Brieskorn homology sphere $\Sigma(2,4,5)$.
\par
If
$u=
\log
\left(
  \frac{1}{2}+\sqrt{2}+\frac{1}{2}\sqrt{5+4\sqrt{2}}
\right)
=1.265948638\ldots$,
then $3u+v(u)=2\pi\sqrt{-1}$, and so this corresponds to the $(3,1)$-surgery along $E$.
We also have $\dfrac{\CS(E_u)}{2\pi^2}=0.041666666667\ldots=\dfrac{1}{24}$, which coincides with the calculation by SnapPea.
In this case $E_{u}$ is $\Sigma(3,3,4)$.
\par
Note that from Figure~\ref{fig:CS} it seems that
\begin{equation*}
  \lim_{u\to\infty}\frac{\CS(E_u)}{2\pi^2}
  =
  \frac{1}{12},
\end{equation*}
which would be the Chern--Simons invariant of $E_{\infty}$ that corresponds to the $(4,1)$-surgery of the figure-eight knot.
\end{rem}
%%%%%%%%%%%%%%%%%%%%%%%%%%%%%%%%%%%%%%%%%%%%%%%%%%%%%%%%%%%%%
\par\medskip\noindent
(ii). When $\theta$ is real and $|\theta|>\arccosh(3/2)$.
\par
In this case $\tilde{H}(u-2\pi\sqrt{-1})$ is real and so $\Im{v(u)}=-2\pi$.
Therefore formulas for $\Vol(E_u)$ and $\CS(E_u)$ similar to the case where $u$ is real and $|u|>\arccosh(3/2)$ hold.
\begin{rem}
Note that $\Re{H(u)}=\Re{\tilde{H}(u-2\pi\sqrt{-1})}$ but $\Im{H(u)}$ and $\Im{\tilde{H}(u-2\pi\sqrt{-1})}$ are different.
In fact $\Im{\tilde{H}(u-2\pi\sqrt{-1})}=0$ but $\Im{H(u)}=2\pi{u}$ for $-\arccosh(3/2)<u<\arccosh(3/2)$.
This means that if $|\Re{u}|\ge\arccosh(3/2)$, then the cases where $\Im{u}=0$ and $\Im{u}=-2\pi\sqrt{-1}$ give different limits but the same $\Vol$ and $\CS$.
As observed in Example~\ref{ex:real_u_volume}, this would correspond to cone-manifold whose underlying manifold is the $0$-surgery along $E$.
\end{rem}
%%%%%%%%%%%%%%%%%%%%%%%%%%%%%%%%%%%%%%%%%%%%%%%%%%%%%%%%%%%%%%%%%%%%%%%%%%%%%
\section{Example~2 -- Torus knots}
In this section we study torus knots.
\par
In particular we will explicitly give representations of the trefoil knot $T:=T(2,3)$ (Figure~\ref{fig:trefoil}) and the cinquefoil knot $C:=T(2,5)$ (Figure~\ref{fig:cinquefoil_pi1_relation}), and consider their relations to the asymptotic behaviors of their colored Jones polynomials.
%%%%%%%%%%%%%%%%%%%%%%%%%%%%%%%%%%%%%%%%%%%%%%%%%%%%%%%%%%%%%%%%%%%%%%%%%%%%%
\subsection{Representations of the fundamental group of the trefoil knot}
\label{subsec:rep_trefoil}
\begin{figure}[h]
\includegraphics[scale=0.5]{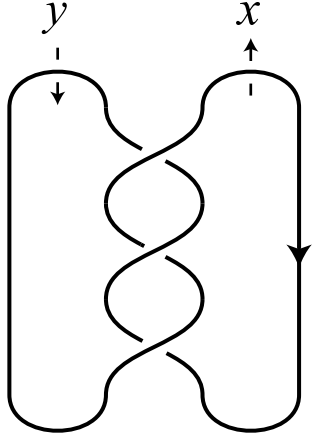}
\caption{The two generators of $\pi_1(S^3\setminus{T})$.}
\label{fig:trefoil_pi1}
\end{figure}

\begin{figure}[h]
\includegraphics[scale=0.5]{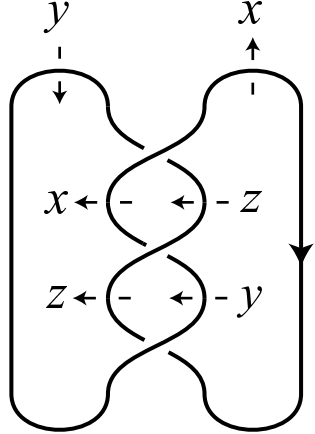}
\caption{The other element of $\pi_1(S^3\setminus{T})$.}
\label{fig:trefoil_pi1_relation}
\end{figure}

We choose the generators $x$ and $y$ for $\pi_1(S^3\setminus{T})$ as in Figure~\ref{fig:trefoil_pi1}.
Note that we take these generators so that their linking numbers with the knot are one.
%%%%%%%%%%%%%
%%%%%%%%%%%%
Then the element $z$ (see Figure~\ref{fig:trefoil_pi1_relation}) can be expressed in terms of $x$ and $y$ as follows:
\begin{equation}\label{eq:trefoil_z}
  z=xyx^{-1}.
\end{equation}
%%%%%%%%%%%%
%%%%%%%%%%%%
From the second crossing and from the third, we have the following relations.
\begin{align*}
  yzx^{-1}z^{-1}&=1,
  \\
  yzy^{-1}z^{-1}&=1.
\end{align*}
Using \eqref{eq:trefoil_z}, these relation are equivalent to the relation
\begin{equation}
  yxy=xyx
\end{equation}
Therefore $\pi_1(S^3\setminus{T})$ has the following presentation.
\begin{equation*}
  \pi_1(S^3\setminus{T})
  =
  \langle
    x,y
    \mid
    \omega x=y\omega
  \rangle
\end{equation*}
with $\omega:=xy$.
\par
Let $\rho$ be a representation of $\pi_1(S^3\setminus{T})$ at $SL(2;\C)$.
By \cite{Riley:QUAJM31984}, we can assume that $\rho$ sends $x$ and $y$ to
\begin{equation*}
  \begin{pmatrix}
    m^{1/2} & 1 \\
    0       & m^{-1/2}
  \end{pmatrix}
  \quad\text{and}\quad
  \begin{pmatrix}
    m^{1/2} & 0 \\
    -d      & m^{-1/2}
  \end{pmatrix},
\end{equation*}
respectively.
Since
\begin{equation*}
  \rho(\omega)
  =
  \begin{pmatrix}
    m-d & m^{-1/2} \\
    -dm^{-1/2} & m^{-1}
  \end{pmatrix},
\end{equation*}
we have
\begin{align*}
  \rho(y\omega)
  &=
  \begin{pmatrix}
    m^{1/2}(m-d) & 1 \\
    -d(m+m^{-1}-d) & m^{-3/2}-dm^{-1/2}
  \end{pmatrix}
  \\
  \intertext{and}
  \rho(\omega x)
  &=
  \begin{pmatrix}
    m^{1/2}(m-d) & m+m^{-1}-d \\
    -d & m^{-3/2}-dm^{-1/2}
  \end{pmatrix}.
\end{align*}
Therefore $d$ should equal $m+m^{-1}-1$ and $m$ uniquely defines a representation.
\par
Putting $\xi:=\tr\bigl(\rho(x)\bigr)$ and $\eta:=\tr\bigl(\rho(xy)\bigr)$, we have
\begin{align*}
  \xi&=m^{1/2}+m^{-1/2} \\
  \eta&=m+m^{-1}-d
\end{align*}
and
\begin{equation*}
  F(\xi,\eta)
  =
  P_{\omega}(\xi,\eta)
  =
  \eta-1.
\end{equation*}
This confirms Theorem~\ref{thm:Le}.
\par
Now reading off from the top-right, the longitude of $\lambda$ is
\begin{equation*}
  \lambda
  =
  yxzx^{-3}
  =
  yx^2yx^{-4}.
\end{equation*}
(Note that we add $x^{-3}$ so that the longitude $\lambda$ has linking number zero with the knot.)
So its image by $\rho$ is
\begin{equation*}
  \rho(\lambda)
  =
  \begin{pmatrix}
    -m^{-3} & \dfrac{m^3-m^{-3}}{m^{1/2}-m^{-1/2}}\\[3mm]
       0    & -m^3
  \end{pmatrix}.
\end{equation*}
\par
We also put $g:=yzx=yxy$ and $h:=xy$.
Then we see that $g^2=h^3$.
So we have another presentation of $\pi_1(S^3\setminus{T})$:
\begin{equation*}
  \langle
    g,h\mid g^2=h^3
  \rangle.
\end{equation*}
We also compute
\begin{align*}
  \rho(g)
  &=
  \begin{pmatrix}
    m^{1/2}-m^{-1/2} & 1 \\
    1-m-m^{-1} & -m^{1/2}+m^{-1/2}
  \end{pmatrix}
  \\
  \intertext{and}
  \rho(h)
  &=
  \begin{pmatrix}
    1-m^{-1} & m^{-1/2} \\
    -m^{1/2}+m^{-1/2}-m^{-3/2} & m^{-1}
  \end{pmatrix}.
\end{align*}
Therefore we have $\tr\bigl(\rho(g)\bigr)=0$ and $\tr\bigl(\rho(h)\bigr)=1$, and so this representation belongs to the same component as $\rho_{1,1}$ in Theorem~\ref{thm:Dubois_Kashaev} in the $SL(2;\C)$-character variety.
%%%%%%%%%%%%%%%%%%%%%%%%%%%%%%%%%%%%%%%%%%%%%%%%%%%%%%%%%%%%%%%
\subsection{Representations of the fundamental group of the cinquefoil knot}
\label{subsec:rep_cinquefoil}
Our next example is the cinquefoil knot (or the double overhand knot) $C=T(2,5)$.
%%%%%%%%%%%%
\begin{figure}[h]
\includegraphics[scale=0.5]{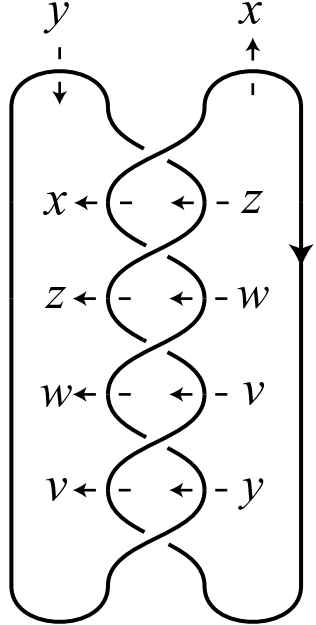}
\caption{Cinquefoil knot and elements of $\pi_1(S^3\setminus{C})$.}
\label{fig:cinquefoil_pi1_relation}
\end{figure}
%%%%%%%%%%%%
The elements $z$, $w$ and $v$ indicated in Figure~\ref{fig:cinquefoil_pi1_relation} are presented in terms of $x$ and $y$:
\begin{align*}
  z&=xyx^{-1},
  \\
  w&=zxz^{-1}=xyxy^{-1}x^{-1},
  \\
  v&=wzw^{-1}=xyxyx^{-1}y^{-1}x^{-1}.
\end{align*}
The last two crossings give the same relation
\begin{equation*}
  yxyxyx^{-1}y^{-1}x^{-1}y^{-1}x^{-1}=1.
\end{equation*}
Therefore we have a presentation of $\pi_1(S^3\setminus{C})$.
\begin{equation*}
  \pi_1(S^3\setminus{C})
  =
  \langle
    x,y
    \mid
    \omega x=y\omega
  \rangle
\end{equation*}
with $\omega=xyxy$.
Then the longitude $\lambda$ is
\begin{equation*}
  \lambda=yxyxy^{-1}xyxyx^{-7}.
\end{equation*}
\par
Let $\rho$ be a non-abelian representation of $\pi_1(S^3\setminus{C})$ at $SL(2,\C)$.
Then we may assume that
\begin{equation*}
  \rho(x)
  =
  \begin{pmatrix}
    m^{1/2} & 1 \\
    0       & m^{-1/2}
  \end{pmatrix}
  \quad\text{and}\quad
  \rho(y)
  =
  \begin{pmatrix}
    m^{1/2} & 0 \\
    -d      & m^{-1/2}
  \end{pmatrix}.
\end{equation*}
Similar calculations show that $d$ and $m$ satisfy the following equation.
\begin{equation}\label{eq:d_m_cinquefoil}
  d^2-(2m+2m^{-1}-1)d+m^2-m+1-m^{-1}+m^{-2}.
\end{equation}
\par
If we put $\xi:=\tr\bigl(\rho(x)\bigr)$ and $\eta:\tr\bigl(xy\bigr)$, then we have
\begin{align*}
  \xi&=m^{1/2}+m^{-1/2},\\
  \eta&=m+m^{-1}-d
\end{align*}
and
\begin{equation*}
  F(\xi,\eta)
  =
  P_{xyxy}(\xi,\eta)-P_{xy}+1
  =
  \eta^2-\eta-1,
\end{equation*}
which coincides \eqref{eq:d_m_cinquefoil} and this confirms Theorem~\ref{thm:Le}.
\par
So we have
\begin{equation*}
  d
  =
  m+m^{-1}-\frac{1\pm\sqrt{5}}{2}
\end{equation*}
and
\begin{equation*}
  \rho(\lambda)
  =
  \begin{pmatrix}
    -m^{-5} & \dfrac{m^5-m^{-5}}{m^{1/2}-m^{-1/2}} \\[3mm]
    0       & -m^5
  \end{pmatrix}
\end{equation*}
for both choices of $d$, where $\rho$ is given as before.
We put $d_{\pm}:=m+m^{-1}-\dfrac{1\pm\sqrt{5}}{2}$ and denote the corresponding representation by $\rho_{\pm}$.
\par
Note that if we put $g:=yvwzx=yxyxy$ and $h:=xy$ then $g^2=h^5$, and so we have another presentation of $\pi_1(S^3\setminus{C})$:
\begin{equation*}
  \langle
    g,h\mid g^2=h^5
  \rangle.
\end{equation*}
We also calculate
\begin{align*}
  \rho_{\pm}(g)
  &=
  \begin{pmatrix}
    \dfrac{1\pm\sqrt{5}}{2}(m^{1/2}-m^{-1/2})
    &
    \dfrac{1\pm\sqrt{5}}{2}
    \\[3mm]
    -\dfrac{1\pm\sqrt{5}}{2}(m-1+m^{-1})+1
    &
    -\dfrac{1\pm\sqrt{5}}{2}(m^{1/2}-m^{-1/2})
  \end{pmatrix},
  \\[5mm]
  \rho_{\pm}(h)
  &=
  \begin{pmatrix}
    \dfrac{1\pm\sqrt{5}}{2}-m^{-1} & m^{-1/2}
    \\
    -m^{-3/2}-m^{1/2}+\dfrac{1\pm\sqrt{5}}{2}m^{-1/2} & m^{-1}
  \end{pmatrix}.
\end{align*}
Therefore we have $\tr\bigl(\rho_{\pm}(g)\bigr)=0$, $\tr\bigl(\rho_{+}(h)\bigr)=\cos(\pi/5)$ and $\tr\bigl(\rho_{-}(h)\bigr)=\cos(3\pi/5)$.
This shows that $\rho_{+}$ and $\rho_{-}$ belong to the same component as $\rho_{1,1}$ and $\rho_{1,3}$ respectively in Theorem~\ref{thm:Dubois_Kashaev}.
%%%%%%%%%%%%%%%%%%%%%%%%%%%%%%%%%%%%%%%%%
\subsection{Asymptotic behavior of the colored Jones polynomial}
In this subsection we review the results in \cite{Murakami:INTJM62004} and \cite{Hikami/Murakami:2007}.
\par
%%%%%%%%%%%%%%%%%%%%%%%%%%%
\begin{thm}\label{thm:torus}
Let $\theta$ be a complex number with $|\theta|>2\pi/(ab)$.
\begin{itemize}
\item
If $\Re(\theta)\Im(\theta)>0$, then we have
\begin{equation*}
  \lim_{N\to\infty}J_N\bigl(T(a,b);\exp(\theta/N)\bigr)
  =
  \frac{1}{\Delta\bigl(T(a,b);\exp(\theta)\bigr)}.
\end{equation*}
\item
If $\Re(\theta)\Im(\theta)<0$, then we have
\begin{equation*}
  \lim_{N\to\infty}
  \frac{\log J_N\bigl(T(a,b);\exp(\theta/N)\bigr)}{N}
  =
  \left(
    1-\frac{\pi\sqrt{-1}}{ab\theta}-\frac{ab\theta}{4\pi\sqrt{-1}}
  \right)
  \pi\sqrt{-1}.
\end{equation*}
\end{itemize}
\end{thm}
%%%%%%%%%%%%%%%%%%%%%%%%%%%
\begin{rem}
In \cite{Murakami:INTJM62004}, the author only proved the case where $\Im(\theta)>0$.
But taking the complex conjugate we have a similar formula for the other case.
This was pointed out by A.~Gibson.
\end{rem}
Putting $u:=\theta-2\pi\sqrt{-1}$, we have
%%%%%%%%%%%%%%%%%%%%%%%%%%%
\begin{cor}\label{cor:torus}
Let $u$ be a complex number with $|u+2\pi\sqrt{-1}|>2\pi/(ab)$.
\begin{itemize}
\item
If $(\Re{u})(\Im{u}+2\pi)>0$, then we have
\begin{equation*}
  \lim_{N\to\infty}J_N\bigl(T(a,b);(u+2\pi\sqrt{-1})/N)\bigr)
  =
  \frac{1}{\Delta\bigl(T(a,b);\exp(u)\bigr)}.
\end{equation*}
\item
If $(\Re{u})(\Im{u}+2\pi)<0$, then we have
\begin{multline*}
  (u+2\pi\sqrt{-1})
  \lim_{N\to\infty}
  \frac{\log J_N\bigl(T(a,b);\exp((u+2\pi\sqrt{-1})/N)\bigr)}{N}
  \\
  =
  \frac{-1}{4ab}\left\{ab(u+2\pi\sqrt{-1})-2\pi\sqrt{-1}\right\}^2.
\end{multline*}
\end{itemize}
\end{cor}
%%%%%%%%%%%%%%%%%%%%%%%%%%%
For small $\theta$, Garoufalidis and L{\^e} proved the following.
%%%%%%%%%%%%%%%%%%%%%%%%%%%
\begin{thm}[{\cite[Theorem~1]{Garoufalidis/Le:aMMR}}]\label{thm:GL_aMMR}
For any knot $K$, there exists a neighborhood $U$ of $0$ such that if $\theta\in U$ then
\begin{equation*}
  \lim_{N\to\infty}J_N\bigl(K;\exp(\theta/N)\bigr)
  =
  \frac{1}{\Delta(K;\exp(\theta))}.
\end{equation*}
\end{thm}
%%%%%%%%%%%%%%%%%%%%%%%%%%%
As in the case of the figure-eight knot when $\exp(\theta)$ a zero of the Alexander polynomial, we can prove that the colored Jones polynomial grows polynomially.
%%%%%%%%%%%%%%%%%%%%%%%%%%%%
\begin{thm}[{\cite[Theorem~1.2]{Hikami/Murakami:2007}}]
We have
\begin{equation*}
  J_N\left(T(a,b);\exp\bigl(\pm2\pi\sqrt{-1}/(abN)\bigr)\right)
  \underset{N\to\infty}{\sim}
  e^{\mp\pi\sqrt{-1}/4}
  \frac{\sin(\pi/a)\sin(\pi/b)}{\sqrt{2}\sin(\pi/(ab))}
  N^{1/2}.
\end{equation*}
\end{thm}
%%%%%%%%%%%%%%%%%%%%%%%%%%%%%%%%%%%%%%%%%%%%
\subsection{Relation of the function $H$ to a representation}
\label{subsec:H_rep_torus}
If $|\theta|$ is small (Theorem~\ref{thm:GL_aMMR}), or $|\theta|>2\pi/(ab)$ and $\Re(\theta)\Im(\theta)>0$ (Theorem~\ref{thm:torus}), then we can associate an abelian representation as in Subsection~\ref{subsec:H_rep_fig8} (iii).
\par
If $|\theta|>2\pi/(ab)$ and $\Re(\theta)\Im(\theta)>0$, that is, if $|u+2\pi\sqrt{-1}|>2\pi/(ab)$ and $(\Re{u})(\Im{u}+2\pi)<0$, we put
\begin{equation*}
  H(u):=\frac{-1}{4ab}\left\{ab(u+2\pi\sqrt{-1})-2\pi\sqrt{-1}\right\}^2
\end{equation*}
and consider its relation to representations at $SL(2;\C)$.
\par
Then the functions $v$ and $f$ defined in Conjecture~\ref{conj:new} become
\begin{equation}\label{eq:v_torus}
  v(u)
  =
  -ab(u+2\pi\sqrt{-1}),
\end{equation}
and
\begin{equation}\label{eq:f_torus}
  f(u)
  =
  -2\pi^2
  +\frac{\pi^2}{ab}
  +ab\pi^2
  -\frac{1}{2}abu\pi\sqrt{-1}
  -H(0).
\end{equation}
Here we leave $H(0)$ since we may need some adjustment.
To consider which representation is associated with $u$, we will recall representations of the trefoil knot and the cinquefoil knot described in Subsections~\ref{subsec:rep_trefoil} and \ref{subsec:rep_cinquefoil}.
%%%%%%%%%%%%%
\subsubsection{Relation of the function $H$ to a representation -- trefoil knot}
Put $m:=\exp(u)$ and consider the representation $\rho$ described in Subsection~\ref{subsec:rep_trefoil}.
Then the meridian $\mu:=x$ is sent to the matrix
\begin{equation*}
  \begin{pmatrix}
    e^{u/2} & 1 \\
      0     & e^{-u/2}
  \end{pmatrix},
\end{equation*}
and the longitude $\lambda$ is to
\begin{equation*}
  \begin{pmatrix}
    -e^{-3u} & \dfrac{e^{3u}-e^{-3u}}{e^{u/2}-e^{-u/2}}\\[3mm]
      0      & -e^{3u}
  \end{pmatrix}
  =
  \begin{pmatrix}
    -e^{v(u)/2} & \ast \\
       0     & -e^{-v(u)/2}
  \end{pmatrix}.
\end{equation*}
Therefore the situation here is the same as the case of the figure-eight knot in Subsection~\ref{subsec:H_rep_fig8}.
%%%%%%%%%%%%%
\subsubsection{Relation of the function $H$ to a representation -- cinquefoil knot}
Similarly we put $m:=\exp(u)$ and consider the representations $\rho_{\pm}$ described in Subsection~\ref{subsec:rep_cinquefoil}.
Then the meridian $\mu:=x$ and the longitude $\lambda$ are sent to
\begin{equation*}
  \begin{pmatrix}
    e^{u/2} & 1 \\
      0     & e^{-u/2}
  \end{pmatrix}
\end{equation*}
and
\begin{equation*}
  \begin{pmatrix}
    -e^{v(u)/2} & \ast \\
       0     & -e^{-v(u)/2}
  \end{pmatrix}
\end{equation*}
respectively.
Therefore both $\rho_{\pm}$ are candidates of such a representation so far.
%%%%%%%%%%%%%%%%%%%%%%%%%%%%%%%%%%%%%%%%%%%%
\subsection{Volume and the Chern--Simons invariant}
The Chern--Simons function of a torus knot for non-abelian part of the $SL(2;\C)$-character variety is calculated in \cite[Proposition~7]{Dubois/Kashaev:MATHA2007}.
Let $\rho_{k,l}$ be the representation described in Theorem~\ref{thm:Dubois_Kashaev}.
Then J.~Dubois and Kashaev proved the following formula \cite[Proposition~4]{Dubois/Kashaev:MATHA2007}:
\begin{multline*}
  \cs_{T(a,b)}\left(\left[\rho_{k,l}\right]\right)
  =
  \\
  \left[
    \frac{u}{4\pi\sqrt{-1}},
    \frac{1}{2}-\frac{abu}{4\pi\sqrt{-1}};
    \exp
    \left(2\pi\sqrt{-1}
      \left(
        \frac{(lad+\varepsilon kbc)^2}{4ab}-\frac{u}{8\pi\sqrt{-1}}
      \right)
    \right)
  \right],
\end{multline*}
where $\varepsilon=\pm$ and the right hand side does not depend of the choice of $\varepsilon$.
\par
We want to express the right hand side by using the same basis as in \eqref{eq:cs_f}.
From Subsection~\ref{subsec:V_CS}, this is equivalent to
\begin{equation}\label{eq:cs_torus}
\begin{split}
  &
  \cs_{T(a,b)}\left(\left[\rho_{k,l}\right]\right)
  \\
  =&
  \left[
    \frac{u}{4\pi\sqrt{-1}},
    \frac{1}{2}-\frac{abu}{4\pi\sqrt{-1}}
    -\frac{ab+1}{2};
  \right.
  \\
  &\qquad
  \left.
    \exp
    \left(
      2\pi\sqrt{-1}
      \left(
        \frac{(lad+\varepsilon kbc)^2}{4ab}-\frac{u}{8\pi\sqrt{-1}}
        +
        \left(\frac{ab+1}{2}\right)\frac{u}{4\pi\sqrt{-1}}
      \right)
    \right)
  \right]
  \\
  =&
  \left[
    \frac{u}{4\pi\sqrt{-1}},
    \frac{-ab(u+2\pi\sqrt{-1})}{4\pi\sqrt{-1}};
    \exp
    \left(
      \frac{\pi\sqrt{-1}(lad+\varepsilon kbc)^2}{2ab}+\frac{abu}{4}
    \right)
  \right]
  \\
  =&
  \left[
    \frac{u}{4\pi\sqrt{-1}},
    \frac{v(u)}{4\pi\sqrt{-1}};
    \exp
    \left(
      \frac{\pi\sqrt{-1}(lad+\varepsilon kbc)^2}{2ab}+\frac{abu}{4}
    \right)
  \right].
\end{split}
\end{equation}
from \eqref{eq:v_torus}.
\par
\begin{rem}
A careful reader may notice that we subtract a half integer $(ab+1)/2$ from the second term, which is not allowed in $SL(2;\C)$ theory.
But it is alright since we are calculating the $PSL(2;\C)$ Chern--Simons invariants (see \cite[p.~543]{Kirk/Klassen:COMMP93}).
\end{rem}
Now we will compare this formula with the $f$ function derived from $H$ in \eqref{eq:f_torus}.
%%%%%%%%%%%%%%%%%%%%%%%%%%%%%%%%%%%%%
\subsubsection{Volume and the Chern--Simons invariant of a representation -- trefoil knot}
Putting $m:=\exp(u)$, $a:=2$, $b:=3$, $c:=-1$, $d:=-1$ and $k=l=1$ in \eqref{eq:cs_torus}, we have
\begin{equation*}
  \cs_{T}\bigl([\rho]\bigr)
  =
  \left[
    \frac{u}{4\pi\sqrt{-1}},
    \frac{3(u+2\pi\sqrt{-1})}{2\pi\sqrt{-1}};
    \exp
    \left(
      \frac{\pi\sqrt{-1}}{12}+\frac{3u}{2}
    \right)
  \right].
\end{equation*}
Therefore the corresponding $f$ function introduced in Subsection~\ref{subsec:V_CS}, which we denote by $\tilde{f}$, is
\begin{equation*}
\begin{split}
  \tilde{f}(u)
  &:=
  \frac{2\pi}{\sqrt{-1}}
  \times
  \left(
     \frac{\pi\sqrt{-1}}{12}+\frac{3u}{2}
  \right)
  \\
  &=
  \frac{\pi^2}{6}-3u\pi\sqrt{-1}.
\end{split}
\end{equation*}
On the other hand, the corresponding $f$ function defined by the $H$ function is:
\begin{equation*}
  f(u)
  =
  \frac{\pi^2}{6}
  -3u\pi\sqrt{-1}
  +4\pi^2
  -H(0)
\end{equation*}
from \eqref{eq:f_torus}.
Therefore $\tilde{f}(u)$ and $f(u)$ coincide modulo $\pi^2$ and $H(0)$.
%%%%%%%%%%%%%%%%%%%%%%%%%%%%%%%%%%%%%%%%%%%%%%%%%%%%%%%%%%%%%%%%%%%%%%%%%%%%%
\subsubsection{Volume and the Chern--Simons invariant of a representation -- cinquefoil knot}
Putting $m=\exp(u)$, $a:=2$, $b:=5$, $c:=-1$, $d:=-2$, we have
\begin{align*}
  \cs_{C}\left(\left[\rho_+\right]\right)
  &=
  \left[
    \frac{u}{4\pi\sqrt{-1}},
    \frac{-5(u+2\pi\sqrt{-1)}}{2\pi\sqrt{-1}};
    \exp
    \left(
      \frac{\pi\sqrt{-1}}{20}+\frac{5u}{2}
    \right)
  \right]
  \\
  \intertext{and}
  \cs_{C}\left(\left[\rho_-\right]\right)
  &=
  \left[
    \frac{u}{4\pi\sqrt{-1}},
    \frac{-5(u+2\pi\sqrt{-1)}}{2\pi\sqrt{-1}};
    \exp
    \left(
      \frac{9\pi\sqrt{-1}}{20}+\frac{5u}{2}
    \right)
  \right],
\end{align*}
where we put $k=l=1$ for $\rho_+$ and $k=1$ and $l=3$ for $\rho_-$.
So the corresponding $\tilde{f}$ functions are
\begin{align*}
  \tilde{f}(u)
  &=
  \frac{\pi^2}{10}-5u\pi\sqrt{-1},
  \\
  \intertext{and}
  \tilde{f}(u)
  &=
  \frac{9\pi^2}{10}-5u\pi\sqrt{-1},
\end{align*}
respectively.
On the other hand the $f$ function defined by the $H$ function is:
\begin{equation*}
  \frac{\pi^2}{10}-5u\pi\sqrt{-1}+8\pi^2-H(0).
\end{equation*}
Therefore if we choose $\rho_{+}$, $f$ and $\tilde{f}$ coincide modulo $\pi^2$ and $H(0)$.
%%%%%%%%%%%%%%%%%%%%%%%%%%%%%
\subsubsection{Volume and the Chern--Simons invariant of a representation -- general torus knot}
Now for a general torus knot $T(a,b)$, let us consider a representation $\rho_{1,1}$ parametrized by $(k,l)=(1,1)$.
Then from \eqref{eq:cs_torus} we have
\begin{equation*}
  \cs_{T(a,b)}\left(\left[\rho_{1,1}\right]\right)
  =
  \left[
    \frac{u}{4\pi\sqrt{-1}},
    \frac{-ab(u+2\pi\sqrt{-1})}{4\pi\sqrt{-1}};
    \exp
    \left(
      \frac{\pi\sqrt{-1}}{2ab}+\frac{abu}{4}
    \right)
  \right],
\end{equation*}
since $ad-bc=1$.
So we have
\begin{equation*}
  \tilde{f}(u)
  =
  \frac{\pi^2}{ab}-\frac{1}{2}abu\pi\sqrt{-1}
\end{equation*}
and
\begin{equation*}
  f(u)
  =
  \frac{\pi^2}{ab}
  -\frac{1}{2}abu\pi\sqrt{-1}
  +(ab-2)\pi^2
  -H(0).
\end{equation*}
Therefore for a general torus knot $f$ and $\tilde{f}$ coincide modulo $\pi^2$ and $H(0)$.
%%%%%%%%%%%%%%%%%%%%%%%%%%%%%%%%%%%%%%%%%%%%%%%%%%%%%%%%%%%%%%%%%%%%%%%%%%%%%
\bibliography{mrabbrev,hitoshi}
\bibliographystyle{hamsplain}
\end{document}